\documentclass{article}

\usepackage{epsfig,xspace}
\usepackage{latexsym}
\usepackage[matrix,arrow,curve]{xy}\xyoption{all}

\usepackage[usenames,dvipsnames,svgnames,table]{xcolor}
\usepackage{amssymb,amsmath,latexsym,amsfonts,amsthm,alltt}
\usepackage[usenames,dvipsnames,svgnames,table]{xcolor}
\usepackage[dvips,usenames]{stmaryrd}
\usepackage{url}
\usepackage{graphicx}
\usepackage{float}
\usepackage{amstext}
\usepackage{cite}
\usepackage{hyperref}
\usepackage{MnSymbol}

\newtheorem{defn}{Definition}[section]
\newtheorem{rem}[defn]{Remark}
\newtheorem{thm}[defn]{Theorem}
\newtheorem{lemma}[defn]{Lemma}
\newtheorem{prop}[defn]{Proposition}
\newtheorem{coro}[defn]{Corollary}

\newtheorem{ex}{Example}[section]

\newcommand{\ra}{\rightarrow}
\newcommand{\lra}{\longrightarrow}

\newcommand{\la}{\leftarrow}

\newcommand{\Ra}{\Rightarrow}

\newcommand{\midsp}{\;|\;}

\newcommand{\telos}{\hfill$\Box$}
\newcommand{\type}[1]{{\tt #1}}

\newcommand{\iso}{\backsimeq}

\newcommand{\upv}{\upVdash}
\newcommand{\rperp}{\mbox{${}^{\upv}$}}

\newcommand{\gphi}{{\mathcal  G}(Y)}
\newcommand{\gpsi}{{\mathcal  G}(X)}

\newcommand{\lperp}{{}\rperp}

\newcommand{\filt}{\mbox{\rm Filt}}

\newcommand{\idl}{\mbox{\rm Idl}}

\title{Reconciliation of Approaches to the Semantics of Logics without Distribution}
\author{Chrysafis (Takis) Hartonas\\
University of Thessaly, Greece\\
$\type{hartonas@uth.gr}$\\[6mm]
{\em In memoriam Jon Michael Dunn}
\date{}}

\begin{document}
\maketitle

\begin{abstract}
This article clarifies and indeed completes an approach (initiated by Dunn and this author several years ago and again pursued by the present author over the last three years or so) to the relational semantics of logics that may lack distribution (Dunn's non-distributive gaggles).  The approach uses sorted frames with an incidence relation on sorts (polarities), equipped with additional sorted relations, but, in the spirit of Occam's razor principle, it drops the extra assumptions made in the generalized Kripke frames approach, initiated by Gehrke, that the frames be separated and reduced (RS-frames). We show in this article that, despite rejecting the additional frame restrictions, all the main ideas and results of the RS-frames approach relating to the semantics of non-distributive logics are captured  in this simpler framework. This contributes in unifying the research field, and, in an important sense, it complements and completes Dunn's gaggle theory project for the particular case of logics that may drop distribution.
\end{abstract}

\section{Preliminaries}
\label{prelims section}

\subsection{A Note on Motivation}
\label{motivation section}
The motivation for this article is twofold. First, it aims at complementing Dunn's gaggle-theory project \cite{dunn-ggl,gglbook} by addressing the case of logics that may lack distribution. This has been the topic of recent research both by the present author and by researchers in the RS-frames approach, initiated by Gehrke \cite{mai-gen}. The second motivation relates to clarifying points of convergence and divergence between our approach and that of RS-frames. Dunn himself seems to have stood at the junction of both, as he has contributed, with this author, the lattice representation result \cite{sdl} on which our approach is based, while having also contributed in applying the RS-frames approach to his gaggle-theory project, with Gehrke and Palmigiano \cite{dunn-gehrke}.  Indeed we conclude that apart from dropping the ``RS'' from the ``RS-frames'' approach, which is to say the assumptions that frames (polarities with relations) are {\em Separated and Reduced}, the two approaches have nearly identical objectives and nearly identical techniques, though they have developed separately. In our concluding remarks we point at some issues that seem to indicate that the approach we have taken may be better suited for some purposes (though no doubt the same can be said for the RS-approach, for some other purposes).

\subsection{Normal Lattice Expansions}
\label{normal lat exp section}
In \cite{dunn-ggl}, Dunn introduced the notion of a distributoid, as a distributive lattice with various operations that in each argument place either distribute or co-distribute over either meets or joins, always returning the same type of lattice operation (always a meet, or always a join). To define technically the extended notion, to which we refer as a {\em normal lattice expansion}, let $\{1,\partial\}$ be a 2-element set, $\mathcal{L}^1=\mathcal{L}$ and $\mathcal{L}^\partial=\mathcal{L}^{op}$ (the opposite lattice, order reversed). Extending the J\'{o}nsson-Tarski terminology \cite{jt1}, a function $f:\mathcal{L}_1\times\cdots\times\mathcal{L}_n\lra\mathcal{L}_{n+1}$ is {\em additive}, or a {\em normal operator}, if it distributes over finite joins of the lattice $\mathcal{L}_i$, for each $i=1,\ldots n$, delivering a join in $\mathcal{L}_{n+1}$.
\begin{defn}\rm
\label{normal lattice op defn}
An $n$-ary operation $f$ on a bounded lattice $\mathcal L$ is {\em a normal lattice operator of distribution type  $\delta(f)=(i_1,\ldots,i_n;i_{n+1})\in\{1,\partial\}^{n+1}$}  if it is a normal additive function  $f:{\mathcal L}^{i_1}\times\cdots\times{\mathcal L}^{i_n}\lra{\mathcal L}^{i_{n+1}}$ (distributing over finite joins in each argument place), where  each $i_j$, for  $j=1,\ldots,n+1$,   is in the set $\{1,\partial\}$, i.e. ${\mathcal L}^{i_j}$ is either $\mathcal L$, or ${\mathcal L}^\partial$.
\end{defn}

\begin{ex}\rm
A modal normal diamond operator $\Diamond$ is a normal lattice operator of distribution type $\delta(\Diamond)=(1;1)$, i.e. $\Diamond:\mathcal{L}\lra\mathcal{L}$, distributing over finite joins of $\mathcal{L}$. A normal box operator $\Box$ is also a normal lattice operator in the sense of Definition \ref{normal lattice op defn}, of distribution type $\delta(\Box)=(\partial;\partial)$, i.e. $\Box:\mathcal{L}^\partial\lra\mathcal{L}^\partial$ distributes over finite joins of $\mathcal{L}^\partial$, which are then just meets of $\mathcal{L}$.

  An {\bf FL}$_{ew}$-algebra (also referred to as a full {\bf BCK}-algebra, or a commutative integral residuated lattice) $\mathcal{A}=(L,\wedge,\vee,0,1,\circ,\ra)$ is a normal lattice expansion, where $\delta(\circ)=(1,1;1)$, $\delta(\ra)=(1,\partial;\partial)$, i.e. $\circ:\mathcal{L}\times\mathcal{L}\lra\mathcal{L}$ and $\ra\;:\mathcal{L}\times\mathcal{L}^\partial\lra\mathcal{L}^\partial$ are both normal lattice operators with the familiar distribution properties.

  Dropping exchange, $\circ$ has two residuals $\la, \ra$, one in each argument place, where $\delta(\la)=(\partial,1;\partial)$, i.e. $\la\;:\mathcal{L}^\partial\times\mathcal{L}\lra\mathcal{L}^\partial$.

  De Morgan Negation $\neg$ is a normal lattice operator and it has both the distribution type $\delta_1(\neg)=(1;\partial)$ and $\delta_2(\neg)=(\partial;1)$, as it switches both joins to meets and meets to joins.

  The Grishin operators \cite{grishin} $\leftharpoondown,\star,\rightharpoondown$, satisfying the familiar co-residuation conditions $a\geq c\leftharpoondown b$ iff $a\star b\geq c$ iff $b\geq a\rightharpoondown c$ have the respective distribution properties, which are exactly captured by assigning to them the distribution types $\delta(\star)=(\partial,\partial;\partial)$ ($\star$ behaves like a binary box operator), $\delta(\leftharpoondown)=(1,\partial;1)$ and $\delta(\rightharpoondown)=(\partial,1;1)$.
\end{ex}
Dunn's distributoid's \cite{dunn-ggl} are the special case of a normal lattice expansion where the underlying lattice is distributive. BAO's (Boolean Algebras with Operators) \cite{jt1,jt2} are the special case where the underlying lattice is a Boolean algebra and all normal operators distribute over finite joins of the Boolean algebra, i.e. they are all of distribution types of the form $(1,\ldots,1;1)$. For BAO's there is no need to consider operators of other distribution types, as they can be obtained by composition of operators with Boolean complementation. For example, in studying residuated Boolean algebras \cite{residBA}, J\'{o}nsson and Tsinakis introduce a notion of {\em conjugate operators} and they show that intensional implications (division operations) $\backslash,/$ (the residuals of the product operator $\circ$)  are interdefinable with the conjugates (at each argument place) $\lhd,\rhd$ of $\circ$, i.e. $a\backslash b=(a\rhd b^-)^-$ and $a\rhd b=(a\backslash b^-)^-$ (and similarly for $/$ and $\lhd$, see \cite{residBA} for details). Note that $\backslash, /$ are not operators, whereas $\lhd,\rhd$ are.

The relational representation of BAO's in \cite{jt1}, extending Stone's representation \cite{stone1} of Boolean algebras using the space of ultrafilters of the algebra, forms the technical basis of the subsequently introduced by Kripke 
possible worlds semantics, with its well-known impact on the development of normal modal logics.
Dunn's approach and objective in \cite{dunn-ggl} has been to achieve the same unified semantic treatment for the logics of distributive lattices with various quasioperators, now based on the Priestley representation \cite{hilary} of distributive lattices in ordered Stone spaces (simplifying Stone's original representation \cite{stone2} of distributive lattices), using the space of prime filters, and abstracting over various specific results in the semantics of distributive, non-classical logics, notably Relevance Logics.

For non-distributive lattices, Urquhart pioneered a representation theorem \cite{urq}, using the space of maximally disjoint filter-ideal pairs. Over the years, Urquhart's representation has proven notoriously difficult to work with, though some authors, including Dunn himself (with Allwein) \cite{dunnLL}, as well as D\"{u}ntsch, Orlowska, Radzikowska and Vakarelov \cite{vac-et-al} have based a semantic treatment of specific systems on it. Hartonas and Dunn \cite{iulg}, published in (1997) \cite{sdl}, provide a lattice representation and duality result based on the representation of semilattices and of Galois connections and abstracting over Goldblatt’s \cite{goldb} representation of ortholattices, replacing orthocomplementation with the trivial Galois connection (the identity map $\imath:\mathcal{L}\lra\left(\mathcal{L}^\partial\right)^\partial$). Hartonas \cite{dloa} presents another lattice representation, extended to include various lattice expansions. Both \cite{sdl,dloa} form the background of a representation and duality result for normal lattice expansions \cite{sdl-exp}, extending the representation of \cite{sdl}.

The bulk of Dunn's work on gaggles predates the extension of the theory of canonical extensions to bounded lattices advanced by Gehrke and Harding in \cite{mai-harding}, extending the J\'{o}nsson-Tarski results for perfect extensions of Boolean algebras \cite{jt1} and the Gehrke-J\'{o}nsson following extension to distributive lattice expansions \cite{mai-jons}. Subsequently, Gehrke \cite{mai-gen} proposed generalized Kripke frames (RS-frames), based on Hartung's lattice representation \cite{hartung}, as a suitable framework for the relational semantics of logics lacking distribution \cite{mai-grishin}, including full Linear Logic \cite{COUMANS201450}. The RS-frames approach to the semantics of logics without distribution was further developed by Palmigiano and co-workers, notably Conradie \cite{conradie2018goldblattthomason,palmigiano-categories,conradie2021nondistributive}.

\section{Polarities with Relations}
\label{polarities section}
\subsection{Definitions, Notational Conventions and Basic Facts}
Let $\{1,\partial\}$ be a set of sorts and $(Z_1,Z_\partial)$ a sorted set.
Base sorted frames $\mathfrak{F}=(X,\upv,Y)=(Z_1,\upv,Z_\partial)$ are triples consisting of nonempty sets $Z_1=X,Z_\partial=Y$ and a binary relation $\upv\;\subseteq X\times Y$, which will be referred to as the {\em Galois relation} of the frame. It generates a Galois connection $(\;)\rperp:\powerset(X)\leftrightarrows\powerset(Y)^\partial:\lperp(\;)$ ($V\subseteq U\rperp$ iff $U\subseteq\lperp V$), defined by
\begin{eqnarray*}
U\rperp&=\{y\in Y\midsp\forall x\in U\; x\upv y\} &= \{y\in Y\midsp U\upv y\}\\
\lperp V&=\{x\in X\midsp \forall y\in V\;x\upv y\}&=\{x\in X\midsp x\upv V\}
\end{eqnarray*}

Triples $(X,R,Y), R\subseteq X\times Y$, where $R$ is treated as the Galois relation of the frame, are variously referred to in the literature as {\em polarities}, after Birkhoff \cite{birkhoff}, as {\em formal contexts}, in the Formal Concept Analysis (FCA) tradition \cite{wille2}, or as {\em object-attribute (categorization, classification, or information) systems} \cite{ewa-rfca,rough-vakarelov}, or as {\em generalized Kripke frames} \cite{mai-gen}, or as {\em polarity frames} in the bi-approximation semantics of \cite{Suzuki-polarity-frames}.

A subset $A\subseteq X$ will be called {\em stable} if $A={}\rperp(A\rperp)$. Similarly, a subset $B\subseteq Y$ will be called {\em co-stable} if $B=({}\rperp B)\rperp$. Stable and co-stable sets will be referred to as {\em Galois sets}, disambiguating to {\em Galois stable} or {\em Galois co-stable} when needed and as appropriate.

$\gpsi,\gphi$ designate the families (complete lattices) of stable and co-stable sets, respectively. Note that the Galois connection restricts to a duality of the complete lattices of Galois stable and co-stable sets $(\;)\rperp:\gpsi\iso\gphi^\partial:{}\rperp(\;)$.

Preorder relations are induced on each of the sorts, by setting for $x,z\in X$, $x\preceq z$ iff $\{x\}\rperp\subseteq\{z\}\rperp$ and, similarly, for $y,v\in Y$, $y\preceq v$ iff ${}\rperp\{y\}\subseteq{}\rperp\{v\}$.
A (sorted) frame is called {\em separated} if the preorders $\preceq$ (on $X$ and on $Y$) are in fact partial orders $\leq$.
Note that if the frame is separated, then $\Gamma x=\Gamma z$ iff $x=z$ and ${}\rperp\{y\}={}\rperp\{v\}$ iff $y=v$. Thus we can identify $X$ and $Y$ with the corresponding subsets of $\gpsi,\gphi$. Moreover, we can identify $Z=X\cup Y$ with the family of sets $\{\Gamma x\midsp x\in X\}\cup\{{}\rperp\{y\}\midsp y\in Y\}\subseteq\gpsi$. The following result is due to Gehrke (\hskip-0.3mm\cite{mai-gen}, Proposition 2.7, Corollary 2.11).
\begin{prop}[Gehrke \cite{mai-gen}]\rm
\label{from Gehrke}
In a separated frame $(X,\upv,Y)$ the set $Z=X\cup Y$ is partially ordered by the relation $\leqslant$ defined for $x,z\in X$ and $y,v\in Y$ by
\begin{tabbing}
\hskip1cm\= $x\leqslant y$\hskip5mm\= iff\hskip2mm\= $\Gamma x\subseteq{}\rperp\{y\}$ \hskip5mm\= iff\hskip2mm\= $x\upv y$\\
\> $x\leqslant z$ \>iff\> $\Gamma x\subseteq\Gamma z$ \>iff\> $z\leq x$\\
\> $y\leqslant v$ \>iff\> ${}\rperp\{y\}\subseteq{}\rperp\{v\}$ \>iff\> $y\leq v$\\
\> $y\leqslant x$ \>iff\> ${}\rperp\{y\}\subseteq\Gamma x$
\>iff\> $\forall u\in X\forall w\in Y\;(u\upv y\;\wedge\;x\upv w\;\lra\;u\upv w)$
\end{tabbing}
Moreover, $\gpsi$ is the Dedekind-MacNeille completion $\overline{Z}$ of $(Z,\leqslant)$.\telos
\end{prop}
Generalized Kripke frames, introduced by Gehrke \cite{mai-gen}, are separated and {\em reduced} polarities (RS-frames), where the latter is defined by the  conditions
\begin{enumerate}
\item $\forall x\in X\exists y\in Y\;(x\not\leqslant y\;\wedge\;\forall z\in X(z\leqslant x\;\wedge\;z\neq x\lra z\leqslant y))$
\item $\forall y\in Y\exists x\in X\;(x\not\leqslant y\;\wedge\;\forall v\in Y(y\leqslant v\wedge y\neq v\lra x\leqslant v))$
\end{enumerate}
The concept goes back to Wille's Formal Concept Analysis (FCA) framework \cite{wille2} and Hartung's lattice representation theorem \cite{hartung}. Gehrke observes that being reduced  means
that all the elements of $X$ are completely join irreducible in $X$ (equivalently, in $Z$, equivalently in $\overline{Z} = \gpsi$) and, dually, that all the elements of $Y$ are completely meet irreducible in $Y$ (equivalently in $Z$, equivalently
in $\overline{Z} = \gpsi$). This turns the poset $(Z,\leqslant)$ to what is called a {\em perfect} poset \cite{dunn-gehrke} (join-generated by the set $J^\infty(Z)$ of its join irreducibles and meet-generated by the set $M^\infty(Z)$ of its meet irreducibles), which can be represented as the two-sorted frame $(J^\infty(Z),M^\infty(Z),\leqslant)$.

To model additional logical operators, RS-frames are equipped with relations subject to the requirement that all their sections be stable sets. For example, to model the Lambek calculus product operator $\circ$, RS-frames are equipped with a relation $S\subseteq Y\times (X\times X)$, all sections of which are required to be stable, and an operation $\bigotimes$ is generated on stable sets by defining (see \cite{mai-gen})
\begin{eqnarray*}
A\mbox{$\bigotimes$}C&=&{}\rperp\{y\in Y\midsp \forall x\in A\forall z\in C\;ySxz\}\\
&=&\{u\in X\midsp\forall y\in Y(\forall x,z\in X(x\in A\wedge z\in C\lra ySxz)\lra u\upv y)\}
\end{eqnarray*}
In the canonical frame (constructed as the Hartung representation \cite{hartung} of the Lindenbaum-Tarski algebra $\mathcal{L}$ of the logic), the relation $S$ is defined by the condition $ySxz$ iff $\forall a,c\in\mathcal{L}\;(a\in x\;\wedge\;c\in z\lra a\circ c\in y)$.

In our own approach, initiated with Dunn and the lattice representation result of \cite{iulg,sdl} and gradually developed in the last three years or so by this author \cite{sdl-exp,pnsds,redm,discres,discr,gts}, we have preferred to apply Occam's razor principle and reject any property of frames that can be rejected while still allowing for the derivation of the needed results. Therefore we work with polarities that need not be separated, or reduced. Section stability for the additional relations, however, is retained as a requirement\footnote{The usefulness of this property has not been fully acknowledged in some of our previous writings, as we did not make section stability explicit. We fill in for this gap in this article.}.

In the rest of this section we specify the dual objects (polarities with relations) of normal lattice expansions, thereby a class of frames for logics that may lack distribution is described. Representation issues (for completeness arguments) is discussed in the next section and it amounts to an extension of the lattice representation published by this author with Dunn \cite{sdl}. We do not specify any particular logical signature (except for assuming that conjunction, disjunction and logical constants for truth and falsity are present). The logical setting is very much the same as that detailed by Conradie and Palmigiano in \cite{conradie-palmigiano} and the reader is referred to this article for a syntactic description of the logics. We save some space by working here with the algebraic structures corresponding to logics that may lack distribution (LE-logics, in the terminology of \cite{conradie2018goldblattthomason}), i.e. with normal lattice expansions.

\begin{rem}[{\bf Notational Conventions}]\rm
\label{notation rem}
For a sorted relation $R\subseteq\prod_{j=1}^{j=n+1}Z_{i_j}$, where $i_j\in\{1,\partial\}$ for each $j$ (and thus $Z_{i_j}=X$ if $i_j=1$ and $Z_{i_j}=Y$ when $i_j=\partial$), we make the convention to regard it as a relation $R\subseteq Z_{i_{n+1}}\times\prod_{j=1}^{j=n}Z_{i_j}$, we agree to write its sort type as $\sigma(R)=(i_{n+1};i_1\cdots i_n)$ and for a tuple of points of suitable sort we write $uRu_1\cdots u_n$ for $(u,u_1,\ldots,u_n)\in R$. We often display the sort type as a superscript, as in $R^\sigma$. Thus, for example, $R^{\partial 1\partial}$ is a subset of  $Y\times(X\times Y)$. In writing then $yR^{\partial 1 \partial}xv$ it is understood that $x\in X=Z_1$ and $y,v\in Y=Z_\partial$. The sort superscript is understood as part of the name designation of the relation, so that, for example, $R^{111}, R^{\partial\partial 1}$ name two different relations.

We use $\Gamma$ to designate  upper closure  $\Gamma U=\{z\in X\midsp\exists x\in U\;x\preceq z\}$, for $U\subseteq X$, and similarly for $U\subseteq Y$. $U$ is {\em increasing} (an upset) iff $U=\Gamma U$. For a singleton set $\{x\}\subseteq X$ we write $\Gamma x$, rather than $\Gamma(\{x\})$ and similarly for $\{y\}\subseteq Y$.

We typically use the standard FCA \cite{wille2} priming notation for each of the two Galois maps ${}\rperp(\;),(\;)\rperp$. This allows for stating and proving results for each of $\gpsi,\gphi$ without either repeating definitions and proofs, or making constant appeals to duality. Thus for a Galois set $G$, $G'=G\rperp$, if $G\in\gpsi$ ($G$ is a Galois stable set), and otherwise $G'={}\rperp G$, if $G\in\gphi$ ($G$ is a Galois co-stable set).

For an element $u$ in either $X$ or $Y$ and a subset $W$, respectively of $Y$ or $X$, we write $u|W$, under a well-sorting assumption, to stand for either $u\upv W$ (which stands for $u\upv w$, for all $w\in W$), or $W\upv u$ (which stands for $w\upv u$, for all $w\in W$), where well-sorting means that either $u\in X, W\subseteq Y$, or $W\subseteq X$ and $u\in Y$, respectively. Similarly for the notation $u|v$, where $u,v$ are elements of different sort.

We designate $n$-tuples (of sets, or elements) using a vectorial notation, setting  $(G_1,\ldots,G_n)=\vec{G}\in\prod_{j=1}^{j=n}\mathcal{G}(Z_{i_j})$, $\vec{U}\in\prod_{j=1}^{j=n}\powerset(Z_{i_j})$, $\vec{u}\in\prod_{j=1}^{j=n}Z_{i_j}$ (where $i_j\in\{1,\partial\}$).   Most of the time we are interested in some particular argument place $1\leq k\leq n$ and we write $\vec{G}[F]_k$ for the tuple $\vec{G}$ where $G_k=F$ (or $G_k$ is replaced by $F$). Similarly $\vec{u}[x]_k$ is $(u_1,\ldots,u_{k-1},x,u_{k+1},\ldots,u_n)$. For brevity, we write $\vec{u}\preceq\vec{v}$ for the pointwise ordering statements $u_1\preceq v_1,\ldots,u_n\preceq v_n$. We also let $\vec{u}\in\vec{W}$ stand for the conjunction of componentwise membership $u_j\in W_j$, for all $j=1,\ldots,n$.

To refer to sections of relations (the sets obtained by leaving one argument place unfilled) we make use of the notation $\vec{u}[\_]_k$ which stands for the $(n-1)$-tuple $(u_1,\ldots,u_{k-1},[\_]\;,u_{k+1},\ldots,u_n)$ and similarly for tuples of sets, extending the membership convention for tuples to cases such as $\vec{u}[\_]_k\in\vec{F}[\_]_k$ and similarly for ordering relations $\vec{u}[\_]_k\preceq\vec{v}[\_]_k$. We also quantify over tuples (with, or without a hole in them), instead of resorting to an iterated quantification over the elements of the tuple, as for example in $\exists\vec{u}[\_]_k\in\vec{F}[\_]_k\exists v,w\in G\;wR\vec{u}[v]_k$.
\end{rem}

\begin{lemma}\rm
\label{basic facts}
Let $\mathfrak{F}=(X,\upv,Y)$ be a  polarity and $u\in Z=X\cup Y$.
\begin{enumerate}
\item $\upv$ is increasing in each argument place
\item $(\Gamma u)'=\{u\}'$ and $\Gamma u=\{u\}^{\prime\prime}$ is a Galois set
\item Galois sets are increasing, i.e. $u\in G$ implies $\Gamma u\subseteq G$
\item For a Galois set $G$, $G=\bigcup_{u\in G}\Gamma u$
\item For a Galois set $G$, $G=\bigvee_{u\in G}\Gamma u=\bigcap_{v|G}\{v\}'$ .
\item For a Galois set $G$ and any set $W$, $W^{\prime\prime}\subseteq G$ iff $W\subseteq G$.
\end{enumerate}
\end{lemma}
\begin{proof}
  By simple calculation. Proof details are included in \cite{sdl-exp}, Lemma 2.2.  For claim 4, $\bigcup_{u\in G}\Gamma u\subseteq G$ by claim 3 (Galois sets are upsets). In claim 5, given our notational conventions, the claim is that if $G\in\gpsi$, then $G=\bigcap_{G\upv y}{}\rperp\{y\}$ and if $G\in\gphi$, then $G=\bigcap_{x\upv G}\{x\}\rperp$.
\end{proof}

For the purposes of this article, the following definition of closed and open elements suffices.
\begin{defn}[Closed and Open Elements]\rm
The principal upper sets of the form $\Gamma x$, with $x\in X$, will be called {\em closed}, or {\em filter} elements of $\gpsi$, while sets of the form ${}\rperp\{y\}$, with $y\in Y$, will be referred to as {\em open}, or {\em ideal} elements of $\gpsi$. Similarly for $\gphi$. A closed element $\Gamma u$ is {\em clopen} iff there exists an element $v$, with $u|v$, such that $\Gamma u=\{v\}'$.
\end{defn}
 By Lemma \ref{basic facts}, the closed elements of $\gpsi$  join-generate $\gpsi$, while the open elements meet-generate $\gpsi$ (similarly for $\gphi$).

\begin{defn}[Galois Dual Relation]\rm\label{Galois dual relations}
For a relation $R$, of sort type $\sigma$, its {\em Galois dual} relation $R'$ is the relation defined by $uR'\vec{v}$ iff $\forall w\;(wR\vec{v}\lra w|u)$. In other words, $R'\vec{v}=(R\vec{v})'$.
\end{defn}
For example, given a relation $R^{111}$ its Galois dual is the relation $R^{\partial 11}$ where for any $x,z\in X$, $R^{\partial 11}xz=(R^{111}xz)\rperp=\{y\in Y\midsp\forall u\in X\;(uR^{111}xz\lra u\upv y)\}$ and, similarly, for a relation $S^{\partial 1\partial}$ its Galois dual is the relation $S^{11\partial}$ where for any $z\in X, v\in Y$ we have $S^{11\partial}zv={}\rperp(S^{\partial 1\partial}zv)$, i.e. $xS^{11\partial}zv$ holds iff for all $y\in Y$, if $yS^{\partial 1\partial}zv$ obtains, then $x\upv y$.

\begin{defn}[Sections of Relations]\rm
\label{sections defn}
For an $(n+1)$-ary relation $R^\sigma$ and an $n$-tuple $\vec{u}$, $R^\sigma\vec{u}=\{w\midsp wR^\sigma\vec{u}\}$ is the {\em section} of $R^\sigma$ determined by $\vec{u}$. To designate a section of the relation at the $k$-th argument place we let $\vec{u}[\_]_k$ be the tuple with a hole at the $k$-th argument place. Then $wR^\sigma\vec{u}[\_]_k=\{v\midsp wR^\sigma\vec{u}[v]_k\}\subseteq Z_{i_k}$ is the $k$-th section of $R^\sigma$.
\end{defn}

\subsection{Image Operators, Conjugates and Residuals}
\label{image ops section}
If $R^\sigma$ is a relation on a sorted frame $\mathfrak{F}$, of some sort type $\sigma=(i_{n+1};i_1\cdots i_n)$, then as in the unsorted case, $R^\sigma$ (but we shall drop the displayed sort type when clear from context)  generates a (sorted) {\em image operator} $\alpha_R$, defined by \eqref{sorted image ops}, of sort $\sigma(\alpha_R)=(i_1,\ldots,i_n;i_{n+1})$, defined by the obvious generalization of the J\'{o}nsson-Tarski image operators \cite{jt1}.
\begin{eqnarray}\label{sorted image ops}
  \alpha_R(\vec{W})&=\;\{w\in Z_{i_{n+1}}\midsp \exists \vec{w}\;(wR\vec{w}\wedge\bigwedge_{j=1}^{j=n}(w_j\in W_j))\} &=\; \bigcup_{\vec{w}\in\vec{W}}R\vec{w}
\end{eqnarray}
where for each $j$, $W_j\subseteq Z_{i_j}$.

Thus $\alpha_R$ is a normal and completely additive function in each argument place, therefore it is residuated, i.e. for each $k$ there is a set-operator $\beta_R^k$ satisfying the condition:
\begin{equation}\label{residuation condition}
\alpha_R(\vec{W}[V]_k)\subseteq U\;\mbox{ iff }\; V\subseteq\beta_R^k(\vec{W}[U]_k)
\end{equation}
Hence $\beta_R^k(\vec{W}[U]_k)$ is the largest set $V$ s.t. $\alpha_R(\vec{W}[V]_k)\subseteq U$ and it is thereby definable by
\begin{equation}\label{def residual of alphaR}
  \beta_R^k(\vec{W}[U]_k)=\bigcup\{V\midsp \alpha_R(\vec{W}[V]_k)\subseteq U\}
\end{equation}
Let $\overline{\alpha}_R$ be the closure of the restriction of $\alpha_R$ to Galois sets $\vec{F}$,
\begin{equation}\label{closure of restriction}
\overline{\alpha}_R(\vec{F})=(\alpha_R(\vec{F}))^{\prime\prime}=\left(\bigcup_{j=1,\ldots,n}^{w_j\in F_j}R\vec{w}\right)^{\prime\prime}=\bigvee_{\vec{w}\in\vec{F}}(R\vec{w})^{\prime\prime}
\end{equation}
where $F_j\in\mathcal{G}(Z_{i_j})$, for each $j\in\{1,\ldots,n\}$.
The operator $\overline{\alpha}_R$ is  sorted  and its sorting is inherited from the sort type of $R$. For example, if $\sigma(R)=(\partial;11)$,  $\alpha_R:\powerset(X)\times\powerset(X)\lra\powerset(Y)$, hence $\overline{\alpha}_R:\gpsi\times\gpsi\lra\gphi$.
Single sorted operations $\overline{\alpha}^1_R:\gpsi\times\gpsi\lra\gpsi$ and $\overline{\alpha}^\partial_R:\gphi\times\gphi\lra\gphi$ can be then extracted by composing appropriately with the Galois connection: $\overline{\alpha}^1_R(A,C)=(\overline{\alpha}_R(A,C))'$ (where $A,C\in\gpsi$) and, similarly, $\overline{\alpha}^\partial_R(B,D)=\overline{\alpha}_R(B',D')$ (where $B,D\in\gphi$). Similarly for the $n$-ary case.

\begin{defn}[Complex Algebra]\rm
\label{complex algebra defn}
Let $\mathfrak{F}=(X,\upv,Y,R)$ be a polarity with a relation $R$ of some sort $\sigma(R)=(i_{n+1};i_1\cdots i_n)$. The {\em complex algebra of} $\mathfrak{F}$ is the structure $\mathfrak{F}^+=(\gpsi,\overline{\alpha}^1_R)$ and its {\em dual complex algebra} is the structure $\mathfrak{F}^\partial=(\gphi,\overline{\alpha}^\partial_R)$.
\end{defn}

\noindent
Most of the time we work with the dual sorted algebra $(\;)\rperp:\gpsi\iso\gphi^\partial:{}\rperp(\;)$, as it allows for considering sorted operations that distribute over joins in each argument place (which are either joins of $\gpsi$, or of $\gphi$, depending on the sort type of the operation). Single-sorted normal operators are then  extracted in the complex algebra by composition with the Galois maps, as indicated above.

\begin{rem}[{\bf Objective}]\rm The primary objective of the current section is to specify conditions under which the residuation structure $\alpha_R\dashv\beta^k_R$ is preserved under the restriction and closure operation described above so that the sorted operator $\overline{\alpha}_R$ on Galois sets is residuated, hence it distributes over arbitrary joins of Galois sets. The notion of conjugate operators we next introduce is useful in this context.
Conjugates were introduced in \cite{residBA} for residuated Boolean algebras. We generalize here to the sorted case, using the duality provided by the Galois connection rather than by classical complementation.
\end{rem}

\begin{defn}[Conjugates]\rm
\label{conjugates defn}
Let $\alpha$ be an image operator (generated by some relation $R$) of sort type $\sigma(\alpha)=(i_1,\ldots,i_k,\ldots,i_n;i_{n+1})$ and $\overline{\alpha}$ the closure of its restriction to Galois sets in each argument place, as defined above.
A function $\overline{\gamma}^k$  on Galois sets, of sort type
\[
\sigma(\overline{\gamma}^k)=(i_1,\ldots,i_{k-1},\overline{i_{n+1}},i_{k+1},\ldots,i_n;\overline{i_k})
\]
(where $\overline{i_j}=\partial$ if $i_j=1$ and $\overline{i_j}=1$ when $i_j=\partial$) is a {\em conjugate} of $\overline{\alpha}$ at the $k$-th argument place (or a $k$-conjugate) iff the following condition holds
\begin{equation}\label{def conjugates}
\hspace*{-2mm}\overline{\alpha}(\vec{F})\subseteq G\;\mbox{ iff }
\overline{\gamma}^k(\vec{F}[G']_k)\subseteq F_k'
\end{equation}
for all Galois sets $F_j\in\mathcal{G}(Z_{i_j})$ and $G\in\mathcal{G}(Z_{i_{n+1}})$.
\end{defn}

It follows from the definition of a conjugate function that $\overline{\gamma}$ is a $k$-conjugate of $\overline{\alpha}$ iff $\overline{\alpha}$ is one of $\overline{\gamma}$ and we thus call $\overline{\alpha},\overline{\gamma}$ $k$-conjugates. Note that the priming notation for both maps of the duality $(\;)\rperp:\gpsi\iso\gphi^\partial:{}\rperp(\;)$ packs together, in one form, four distinct (due to sorting) cases of conjugacy.

\begin{ex}\rm
\label{example and strategy ex 1}
In the case of a ternary relation $R^{111}$ of the indicated sort type, an image operator $\alpha_R=\bigodot:\powerset(X)\times\powerset(X)\lra\powerset(X)$ is generated.  Designate the closure of its restriction to Galois stable sets by $\bigovert:\gpsi\times\gpsi\lra\gpsi$.
Then $\overline{\alpha}=\bigovert$ is of sort type $\sigma(\bigovert)=(1,1;1)$. If $\overline{\gamma}^2_R=\triangleright:\gpsi\times\gphi\lra\gphi$, with $\sigma(\triangleright)=(1,\partial;\partial)$,  then $\bigovert$, $\triangleright$ are {\em conjugates} iff for any Galois stable sets $A,F,C\in\gpsi$ it holds that $A\bigovert F\subseteq C$ iff $A\triangleright C'\subseteq F'$.

Note that, given an operator $\triangleright:\gpsi\times\gphi\lra\gphi$, if we now define $\Ra\;:\gpsi\times\gpsi\lra\gpsi$ by $A\Ra C=(A\triangleright C')'={}\rperp(A\triangleright C\rperp)$, it is immediate that $\bigovert, \triangleright$ are conjugates iff $\bigovert,\Ra$ are residuated. In other words
\[
A\mbox{$\bigovert$}F\subseteq C\;\mbox{ iff }\;A\triangleright C'\subseteq F'\;\mbox{ iff }\;F\subseteq A\Ra C
\]
\end{ex}

\begin{lemma}\rm\label{conjugate-residual}
The following are equivalent.
\begin{enumerate}
\item[1)] $\overline{\alpha}_R$ distributes over any joins of Galois sets at the $k$-th argument place
\item[2)] $\overline{\alpha}_R$ has a $k$-conjugate $\overline{\gamma}_R^k$ defined on Galois sets by
  \[
  \overline{\gamma}_R^k(\vec{F})=\bigcap\{G\midsp \overline{\alpha}_R(\vec{F}[G']_k)\subseteq F'_k\}
  \]
\item[3)] $\overline{\alpha}_R$ has a $k$-residual $\overline{\beta}_R^k$  defined on Galois sets by
  \[
      \overline{\beta}_R^k(\vec{F}[G]_k)= (\overline{\gamma}_R^k(\vec{F}[G']_k))'=\bigvee\{G^\prime\midsp \overline{\alpha}_R(\vec{F}[G']_k)\subseteq F'_k\}
  \]
\end{enumerate}
\end{lemma}
\begin{proof}
Existence of a $k$-residual is equivalent to distribution over arbitrary joins and the residual is defined by
\[
\overline{\beta}^k_R(\ldots,F_{k-1},H,F_{k+1},\ldots)=\bigvee\{G\midsp\overline{\alpha}_R(\ldots,F_{k-1},G,F_{k+1},\ldots)\subseteq H\}
\]
We show that the distributivity assumption 1) implies that 2) and 3) are equivalent, i.e. that
\[
\overline{\alpha}_R(\vec{F}[G]_k)\subseteq H\;\mbox{ iff }\;
\overline{\gamma}^k_R(\vec{F}[H']_k)\subseteq G'\;\mbox{ iff }\;G\subseteq\overline{\beta}^k_R(\vec{F}[H]_k)
\]

We illustrate the proof for the unary case only, as the other parameters remain idle in the argument.

Assume $\overline{\alpha}_R(G)\subseteq H$ and let $\overline{\gamma}_R(H^\prime)=\bigcap\{E\midsp\overline{\alpha}_R(E')\subseteq H\}$, a Galois set by definition, given that $G,H,E$ are assumed to be Galois sets. Then $G^\prime$ is in the set whose intersection is taken. Hence $\overline{\gamma}_R(H^\prime)\subseteq G^\prime$ follows from the definition of $\overline{\gamma}_R$. It also follows by  definition that $G\subseteq\overline{\beta}_R(H)=(\overline{\gamma}_R(H^\prime))^\prime$.

Assuming $G\subseteq\overline{\beta}_R(H)$ we obtain by definition that $G\subseteq (\overline{\gamma}_R(H^\prime))^\prime$, hence  $G\subseteq\bigvee\{E^\prime\midsp\overline{\alpha}_R(E^\prime)\subseteq H\}$, using the definition of $\overline{\gamma}_R$ and duality. Hence by the distributivity assumption $\overline{\alpha}_R(G)\subseteq\bigvee\{\overline{\alpha}_R(E^\prime)\midsp \overline{\alpha}_R(E^\prime)\subseteq H\}\subseteq H$. This establishes that
$\overline{\alpha}_R(G)\subseteq H$ iff $\overline{\gamma}_R(H^\prime)\subseteq G^\prime$ iff $G\subseteq\overline{\beta}_R(H)$, qed.
\end{proof}

\begin{defn}\rm
We let $\beta^k_{R/}$ be the restriction of $\beta^k_R$ of equation \eqref{def residual of alphaR} to Galois sets, according to its sort type, explicitly defined by \eqref{betakR}:

\begin{equation}\label{betakR}
  \beta^k_{R/}(\vec{E}[G]_k)=\bigcup\{F\in\mathcal{G}(Z_{i_k})\midsp\alpha_R(\vec{E}[F]_k)\subseteq G\}
\end{equation}
\end{defn}

\begin{thm}
\label{beta1R residual}
  \rm
  If $\overline{\alpha}_R$ is residuated in the $k$-th argument place, then $\beta^k_{R/}$ is its residual and $ \beta^k_{R/}(\vec{E}[G]_k)$ is a Galois set, i.e. the union in equation \eqref{betakR} is actually a join in $\mathcal{G}(Z_{i_k})$.
\end{thm}
\begin{proof}
  We illustrate the proof for the unary case only, since the other parameters that may exist remain idle in the argument. In the unary case, $\beta_{R/}(G)=\bigcup\{F\midsp\alpha_R(F)\subseteq G\}$, for Galois sets $F,G$.

  Note first that $\overline{\alpha}_R(F)\subseteq G$ iff $F\subseteq\beta_{R/}(G)$. Left-to-right is obvious by definition and by the fact that for a Galois set $G$ and any set $U$, $U^{\prime\prime}\subseteq G$ iff $U\subseteq G$.  If $F\subseteq\beta_{R/}(G)\subseteq\beta_R(G)$, then by residuation $\alpha_R(F)\subseteq G$. Given that $G$ is a Galois set, it follows $\overline{\alpha}_R(F)\subseteq G$.

 If indeed $\overline{\alpha}_R$ is residuated on Galois sets with a map $\overline{\beta}_R$, then the residual is defined by $\overline{\beta}_R(G)=\bigvee\{F\midsp \overline{\alpha}_R(F)\subseteq G\}=\bigvee\{F\midsp \alpha_R(F)\subseteq G\}$ and this is precisely the closure of $\beta_{R/}(G)=\bigcup\{F\midsp\alpha_R(F)\subseteq G\}$. But in that case we obtain $F\subseteq\overline{\beta}_R(G)$ iff $\overline{\alpha}_R(F)\subseteq G$ iff $\alpha_R(F)\subseteq G$ iff $F\subseteq\beta_{R/}(G)$ and setting $F=\overline{\beta}_R(G)$ it follows that $\overline{\beta}_R(G)\subseteq\beta_{R/}(G)\subseteq\overline{\beta}_R(G)$.
\end{proof}

\begin{lemma}\rm
\label{equiv defn of betakr}
$\beta^k_{R/}$ is equivalently defined by \eqref{beta equiv 1} and by \eqref{beta equiv 2}
\begin{eqnarray}
\beta^k_{R/}(\vec{E}[G]_k)&=&\bigcup\{\Gamma u\in\mathcal{G}(Z_{i_k})\midsp\alpha_R(\vec{E}[\Gamma u]_k)\subseteq G\}
\label{beta equiv 1}\\
\beta^k_{R/}(\vec{E}[G]_k)&=&\{ u\in Z_{i_k}\midsp\alpha_R(\vec{E}[\Gamma u]_k)\subseteq G\}
\label{beta equiv 2}
\end{eqnarray}
\end{lemma}
\begin{proof}
  $\beta^k_{R/}$ is defined by equation \eqref{betakR}, so if $u\in \beta^k_{R/}(\vec{E}[G]_k)$, let $F\in \mathcal{G}(Z_{i_k})$ be such that $u\in F$ and $\alpha_R(\vec{E}[F]_k)\subseteq G$. Then $\Gamma u\subseteq F$ and by monotonicity of $\alpha_R$ we have
  $
  \alpha_R(\vec{E}[\Gamma u]_k)\subseteq\alpha_R(\vec{E}[F]_k)\subseteq G
   $
   and this establishes the left-to-right inclusion for the first identity of the lemma. The converse inclusion is obvious since $\Gamma u$ is a Galois set.

  For the second identity, the inclusion right-to-left is obvious. Now if $u$ is such that $\alpha_R(\vec{E}[\Gamma u]_k)\subseteq G$ and $u\preceq w$, then $\Gamma w\subseteq\Gamma u$ and then by monotonicity of $\alpha_R$ it follows that $\alpha_R(\vec{E}[\Gamma w]_k)\subseteq \alpha_R(\vec{E}[\Gamma u]_k)\subseteq G$.

  This shows that  $\bigcup\{\Gamma u\in\mathcal{G}(Z_{i_k})\midsp\alpha_R(\vec{E}[\Gamma u]_k)\subseteq G\}$ is contained in the set $\{ u\in Z_{i_k}\midsp\alpha_R(\vec{E}[\Gamma u]_k)\subseteq G\}$, and given the first part of the lemma, the second identity obtains as well.
\end{proof}

\begin{defn}[Conjugate Relations]\rm
\label{conjugate relations defn}
Let $\mathfrak{F}=(X,\upv,Y,R,S)$, where $\sigma(R)=(i_{n+1};i_1\cdots i_k\cdots i_n)$, $\sigma(S)=(t_{n+1};t_1\cdots t_k\cdots t_n)$, where $t_{n+1}=\overline{i_k}, t_k=\overline{i_{n+1}}$ and for $j\not\in\{k,n+1\}$, $t_j=i_j$. Let $\alpha_R$ and $\eta_S$ be the generated image operators and $\overline{\alpha}_R,\overline{\eta}_S$ be the closures of their restriction to Galois sets.

The relations $R,S$ will be called {\em $k$-conjugate relations} iff the Galois set operators $\overline{\alpha}_R,\overline{\eta}_S$ are $k$-conjugates (Definition \ref{conjugates defn}), i.e. just in case (given that $G,F'_k$ are Galois sets) $\alpha_R(\vec{F})\subseteq G\;\mbox{ iff }\; \eta_S(\vec{F}[G']_k)\subseteq F^\prime_k$.
\end{defn}

\begin{lemma}\rm
\label{condition for conjugate relations}
Let $\mathfrak{F}=(X,\upv,Y,R,S)$, where $\sigma(R)=(i_{n+1};i_1\cdots i_k\cdots i_n)$, $\sigma(S)=(t_{n+1};t_1\cdots t_k\cdots t_n)$, where $t_{n+1}=\overline{i_k}, t_k=\overline{i_{n+1}}$ and for $j\not\in\{k,n+1\}$, $t_j=i_j$.
Assume that the $k$-th sections of the Galois dual relation $R'$ of $R$ are Galois sets.
Let $T$ be the relation defined, for $w\in Z_{\overline{i_k}}$, by
 $vT\vec{p}[w]_k\;\mbox{ iff }\; w\in (vR'\vec{p}[\_]_k)'$ iff $\forall u\in F_k\;(vR'\vec{p}[u]_k\lra u|w)$.

 If the constraint \eqref{k-conjugate relations constraint} holds in the frame, then $R$ and $S$ are $k$-conjugates.
 \begin{equation}\label{k-conjugate relations constraint}
   \forall v\in Z_{\overline{i_{n+1}}}\forall \vec{p}[\_]_k\in\vec{Z_{i_j}}[\_]_k\forall w\in Z_{\overline{i_k}}\left( vT\vec{p}[w]_k\;\leftrightarrow\;wS\vec{p}[v]_k \right)
 \end{equation}
\end{lemma}
\begin{proof}
We have
\begin{tabbing}
$\alpha_R(\vec{F})\subseteq G$\hskip5mm\= iff\hskip5mm\= $\bigcup_{\vec{p}\in\vec{F}}R\vec{p}\;\subseteq G$ \hskip10mm\= iff\hskip2mm\=$\forall\vec{p}\;(\vec{p}\in\vec{F}\;\lra\; (R\vec{p}\subseteq G))$\\
\hskip5mm\=iff\hskip2mm\= $\forall\vec{p}\;(\vec{p}\in\vec{F}\;\lra\;(G'\subseteq R'\vec{p}))$\\
 \>iff\> $\forall\vec{p}\;(\vec{p}\in\vec{F}\;\lra\;\forall v\in Z_{\overline{i_{n+1}}}(G|v\lra vR'\vec{p}))$\\
\>iff\> $\forall\vec{p}\forall v\in Z_{\overline{i_{n+1}}}\;(\vec{p}[\_]_k\in\vec{F}[\_]_k\wedge p_k\in F_k\wedge G|v\;\lra\; vR'\vec{p}[p_k]_k)$\\
\>iff\>  $\forall\vec{p}\forall v\in Z_{\overline{i_{n+1}}}\;(\vec{p}[\_]_k\in\vec{F}[\_]_k \wedge G|v\;\lra\; (p_k\in F_k\lra vR'\vec{p}[p_k]_k))$\\
\>iff\> $\forall\vec{p}[\_]_k\forall v\in Z_{\overline{i_{n+1}}}\;(\vec{p}[\_]_k\in\vec{F}[\_]_k \wedge G|v\;\lra\; (F_k\subseteq vR'\vec{p}[\_]_k))$\\
\>\> (using the hypothesis that the $k$-th sections of $R'$ are Galois sets)\\
\>iff\> $\forall\vec{p}[\_]_k\forall v\in Z_{\overline{i_{n+1}}}\;(\vec{p}[\_]_k\in\vec{F}[\_]_k \wedge G|v\;\lra\; (\;(vR'\vec{p}[\_]_k)'\subseteq F^\prime_k))$\\
\>iff\> $\forall\vec{p}[\_]_k\forall v\in Z_{\overline{i_{n+1}}}\;\left(\vec{p}[\_]_k\in\vec{F}[\_]_k \wedge G|v\;\lra\; \forall w\in Z_{\overline{i_k}} (\;vT\vec{p}[w]_k\lra F_k|w)\right)$\\
\>iff\> $\forall\vec{p}[\_]_k\forall v\in Z_{\overline{i_{n+1}}}\forall w\in Z_{\overline{i_k}}\; \left(vT\vec{p}[w]_k \wedge \vec{p}[\_]_k\in\vec{F}[\_]_k \wedge G|v\;\lra\; F_k|w\right)$
\end{tabbing}
On the other hand, we have
\begin{tabbing}
$\eta_S(\vec{F}[G']_k)\subseteq F^\prime_k$ \hskip5mm\= iff\hskip5mm\= $\bigcup_{\vec{p}[v]_k\in\vec{F}[G']_k} S\vec{p}[v]_k\;\subseteq F^\prime_k$\\
\hskip5mm\=iff\hskip2mm\= $\forall\vec{p}[\_]_k\forall v\in Z_{\overline{i_{n+1}}}\forall w\in Z_{\overline{i_k}}\; \left(wS\vec{p}[v]_k \wedge \vec{p}[\_]_k\in\vec{F}[\_]_k \wedge G|v\;\lra\; F_k|w\right)$
\end{tabbing}
and thus the claim of the lemma is proved.
\end{proof}

\begin{thm}\rm
\label{distribution from stability thm}
Let $\mathfrak{F}=(X,\upv,Y,R)$ be a frame with an $(n+1)$-ary sorted relation, of some sort  $\sigma(R)=(i_{n+1};i_1\cdots i_n)$. If for any $w\in Z_{\overline{i_{n+1}}}$ and any $(n-1)$-tuple $\vec{p}[\_]_k$ with $p_j\in Z_{i_j}$, for each $j\in\{1,\ldots,n\}\setminus\{k\}$, the sections $wR'\vec{p}[\_]_k$ of the Galois dual relation $R'$ of $R$ are Galois sets, then $\overline{\alpha}_R$ distributes at the $k$-th argument place over arbitrary joins in $\mathcal{G}(Z_{i_k})$.
\end{thm}
\begin{proof}
  Define the relation $T$ from $R$ as in the statement of Lemma \ref{condition for conjugate relations},
  \[
  vT\vec{p}[w]_k\;\mbox{ iff }\; w\in (vR'\vec{p}[\_]_k)'
  \]
  Then use equation \eqref{k-conjugate relations constraint}, repeated below, as a definition for a relation $S$
   \begin{eqnarray*}
   \forall v\in Z_{\overline{i_{n+1}}}\forall \vec{p}[\_]_k\in\vec{Z_{i_j}}[\_]_k\forall w\in Z_{\overline{i_k}}\left( vT\vec{p}[w]_k\;\leftrightarrow\;wS\vec{p}[v]_k \right)
 \end{eqnarray*}
 Note that the sort type of $S$, as defined, is $\sigma(S)=(t_{n+1};t_1\cdots t_k\cdots t_n)$, where $t_{n+1}=\overline{i_k}, t_k=\overline{i_{n+1}}$ and for $j\not\in\{k,n+1\}$, $t_j=i_j$.
 By the proof of Lemma \ref{condition for conjugate relations}, the relations $R$ and $S$ are $k$-conjugates. Consequently, by Lemma \ref{conjugate-residual}, $\overline{\alpha}_R$ distributes at the $k$-th argument place over arbitrary joins in $\mathcal{G}(Z_{i_k})$ and it has a $k$-residual which, by Theorem \ref{beta1R residual} is precisely the restriction to Galois sets $\beta^k_{R/}$ (defined by equation \eqref{betakR}, equivalently by Lemma \ref{equiv defn of betakr}) of the $k$-residual $\beta^k_R$ of the image operator $\alpha_R$.
\end{proof}

By composition with the Galois connection, single-sorted operators $\overline{\alpha}^1_R,\overline{\alpha}^\partial_R$ can be obtained on $\gpsi$ and $\gphi$, respectively. Given that the Galois connection is a duality between Galois stable and co-stable sets, completely normal lattice operators (dual to each other) are obtained on $\gpsi$ and $\gphi$, respectively.
Therefore we have proven the following result.
\begin{coro}\rm
\label{from polarities to lattice expansions}
Let $\mathfrak{F}=(X,\upv,Y,(R_p)_{p\in P})$ be a polarity with relations indexed in some set $P$, of varying sort types $\sigma_p\;(p\in P)$ and such that every section of the Galois dual relations $R'_p\;(p\in P)$ is a Galois set. Then the dual complex algebra $\mathfrak{F}^+$ of $\mathfrak{F}$ is a normal lattice expansion where a relation of sort type $(i_{n+1};i_1\cdots i_n)$ determines a completely normal lattice operator of distribution type $(i_1,\ldots,i_n;i_{n+1})$.\telos
\end{coro}

\section{Representation of Normal Lattice Expansions}
\label{representation section}
A bounded lattice expansion is a structure $\mathcal{L}=(L,\leq,\wedge,\vee,0,1,\mathcal{F}_1,\mathcal{F}_\partial)$, where $\mathcal{F}_1$ consists of normal lattice operators $f$ of distribution type $\delta(f)=(i_1,\ldots,i_n;1)$ (i.e. of output type 1), while $\mathcal{F}_\partial$  consists of normal lattice operators $h$ of distribution type $\delta(h)=(t_1,\ldots,t_n;\partial)$ (i.e. of output type $\partial$). For representation purposes, nothing depends on the size of the operator families $\mathcal{F}_1$ and $\mathcal{F}_\partial$ and we may as well assume that they contain a single member, say $\mathcal{F}_1=\{f\}$ and $\mathcal{F}_\partial=\{h\}$. In addition, nothing depends on the arity of the operators, so we may assume they are both $n$-ary.

\subsection{Canonical Frame Construction}
\label{canonical frame section}
The canonical frame is constructed as follows, based on \cite{iulg,sdl,dloa,sdl-exp}.\\

First, the base polarity $\mathfrak{F}=(\filt(\mathcal{L}),\upv,\idl(\mathcal{L}))$ consists of the sets $X=\filt(\mathcal{L})$ of filters and $Y=\idl(\mathcal{L})$ of ideals of the lattice and the relation $\upv\;\subseteq\filt(\mathcal{L})\times\idl(\mathcal{L})$ is defined by $x\upv y$ iff $x\cap y\neq\emptyset$, while the representation map $\zeta_1$ sends a lattice element $a\in L$ to the set of filters that contain it, $\zeta_1(a)=\{x\in X\midsp a\in x\}=\{x\in X\midsp x_a\subseteq x\}=\Gamma x_a$. Similarly, a co-represenation map $\zeta_\partial$ is defined by $\zeta_\partial(a)=\{y\in Y\midsp a\in y\}=\{y\in Y\midsp y_a\subseteq y\}=\Gamma y_a$. It is easily seen that $(\zeta_1(a))'=\zeta_\partial(a)$ and, similarly, $(\zeta_\partial(a))'=\zeta_1(a)$. The images of $\zeta_1,\zeta_\partial$ are precisely the families (sublattices of $\gpsi,\gphi$, respectively) of clopen elements of $\gpsi,\gphi$, since clearly $\Gamma x_a={}\rperp\{y_a\}$ and $\Gamma y_a=\{x_a\}\rperp$. For further details the reader is referred to \cite{iulg,sdl}.

Second, for each normal lattice operator a relation is defined, such that if $\delta=(i_1,\ldots,i_n;i_{n+1})$ is the distribution type of the operator, then $\sigma=(i_{n+1};i_1\cdots i_n)$ is the sort type of the relation. Without loss of generality, we have restricted to the families of operators $\mathcal{F}_1=\{f\}$ and $\mathcal{F}_\partial=\{h\}$, so that we shall define two corresponding relations $R,S$ of respective sort types $\sigma(R)=(1;i_1\cdots i_n)$ and $\sigma(S)=(\partial;t_1\cdots t_n)$, where for each $j$, $i_j$ and $t_j$ are in $\{1,\partial\}$. In other words
\[
R\;\subseteq\; X\times\prod_{j=1}^{j=n}Z_{i_j}\hskip2cm S\;\subseteq\; Y\times \prod_{j=1}^{j=n}Z_{t_j}
\]
To define the relations, we use the point operators introduced in \cite{dloa} (see also \cite{sdl-exp}). In the generic case we examine, we need to define two sorted operators
\begin{eqnarray*}
\widehat{f}&:\prod_{j=1}^{j=n}Z_{i_j}\lra Z_1\hskip1cm
\widehat{h}&: \prod_{j=1}^{j=n}Z_{t_j}\lra Z_\partial\hskip1cm(\mbox{recall that }Z_1=X, Z_\partial=Y)
\end{eqnarray*}
Assuming for the moment that the point operators have been defined, the canonical relations $R,S$ are defined by
\begin{eqnarray}
xR\vec{u} & \mbox{ iff }& \widehat{f}(\vec{u})\subseteq x\;\; (\mbox{for }\; x\in X\;\mbox{ and }\;\vec{u}\in \prod_{j=1}^{j=n}Z_{i_j})\nonumber\\
yS\vec{v} &\mbox{ iff }& \widehat{h}(\vec{v})\subseteq y\;\; (\mbox{for }\; y\in Y\;\mbox{ and }\;\vec{v}\in \prod_{j=1}^{j=n}Z_{t_j})\label{canonical relations defn}
\end{eqnarray}
Returning to the point operators and letting $x_e,y_e$ be the principal filter and principal ideal, respectively, generated by a lattice element $e$, these are uniformly defined as follows, for $\vec{u}\in \prod_{j=1}^{j=n}Z_{i_j}$ and $\vec{v}\in \prod_{j=1}^{j=n}Z_{t_j})$
\begin{eqnarray}
  \widehat{f}(u_1,\ldots,u_n) &=& \bigvee\{x_{f(a_1,\ldots,a_n)}\midsp \bigwedge_j\;(a_j\in u_j)\}=\bigvee\{x_{f(\vec{a})}\midsp\vec{a}\in\vec{u}\}\nonumber \\
  \widehat{h}(v_1,\ldots,v_n) &=& \bigvee\{y_{h(a_1,\ldots,a_n)}\midsp \bigwedge_j\;(a_j\in v_j)\} =\bigvee\{y_{h(\vec{a})}\midsp\vec{a}\in\vec{v}\}\label{canonical point operators defn}
\end{eqnarray}
In other words, $\widehat{f}(\vec{u})=\langle\{f(\vec{a})\midsp \vec{a}\in\vec{u}\}\rangle$ is the filter generated by the set $\{f(\vec{a})\midsp \vec{a}\in\vec{u}\}$ and similarly $\widehat{h}(\vec{v})$ is the ideal generated by the set $\{h(\vec{a})\midsp \vec{a}\in\vec{v}\}$.

\begin{ex}[FL$_{ew}$]\rm
\label{example and strategy ex 3}
We consider as an example the case of associative, commutative, integral residuated lattices $\mathcal{L}=(L,\leq,\wedge,\vee,0,1,\circ,\ra)$, the algebraic models of {\bf FL}$_{ew}$ (the associative full Lambek calculus with exchange and weakening), also referred to in the literature as full {\bf BCK}. By residuation of $\circ,\ra$, the distribution types of the operators are $\delta(\circ)=(1,1;1)$ and $\delta(\ra)=(1,\partial;\partial)$. Let $(\filt(\mathcal{L}),\upv,\idl(\mathcal{L}))$ be the canonical frame of the bounded lattice $(L,\leq,\wedge,\vee,0,1)$. Designate the corresponding canonical point operators by $\overt$ and $\leadsto$, respectively. They are defined by \eqref{canonical point   operators defn}
\begin{eqnarray*}
x\overt z &=& \bigvee\{x_{a\circ c}\midsp a\in x\;\wedge\; c\in z\}\in\filt(\mathcal{L})\hskip1cm(x,z\in\filt(\mathcal{L}))\\
x\leadsto v &=&\bigvee\{y_{a\ra c}\midsp a\in x\;\wedge\; c\in v\}\in\idl(\mathcal{L})\hskip9.3mm (x\in\filt(\mathcal{L}),v\in\idl(\mathcal{L}))
\end{eqnarray*}
where recall that we write $x_e,y_e$ for the principal filter and ideal, respectively, generated by the lattice element $e$, so that $x\overt z\in \filt(\mathcal{L})$, while $(x\leadsto v)\in \idl(\mathcal{L})$.

The relations $R^{111},S^{\partial 1\partial}$ are then defined by
\[
uR^{111}xz\;\mbox{ iff }\; x\overt z\subseteq u\hskip2cm yS^{\partial 1\partial}xv\;\mbox{ iff }\;(x\leadsto v)\subseteq y
\]
of sort types $\sigma(R)=(1;11)$ and $\sigma(S)=(\partial;1\partial)$. The canonical {\bf FL}$_{ew}$-frame is therefore the structure $\mathfrak{F}=(\filt(\mathcal{L}),\upv,\idl(\mathcal{L}),R^{111}, S^{\partial 1\partial})$.
\end{ex}

\subsection{Properties of the Canonical Frame}
\begin{lemma}\rm
\label{elementary props in canonical}
The following hold for the canonical frame.
\begin{enumerate}
\item The frame is separated
\item For $\vec{u}\in \prod_{j=1}^{j=n}Z_{i_j}$ and $\vec{v}\in \prod_{j=1}^{j=n}Z_{t_j}$ the sections $R\vec{u}$ and $S\vec{v}$ are closed elements of $\gpsi$ and $\gphi$, respectively
\item For $x\in X, y\in Y$, the $n$-ary relations $xR, yS$ are decreasing in every argument place
\end{enumerate}
\begin{proof}
  \mbox{}
For 1),
just note that the ordering $\preceq$ is set-theoretic inclusion (of filters, and of ideals, respectively), hence separation of the frame is immediate.

For 2),  by the definition of the relations,  $R\vec{u}=\{x\midsp \widehat{f}(\vec{u})\subseteq x\}=\Gamma(\widehat{f}(\vec{u}))$ is a closed element of $\gpsi$ and similarly for $S\vec{v}$.

For 3),  if $w\subseteq u_k$, then $\{x_{f(a_1,\ldots,a_n)}\midsp a_k\in w\wedge \bigwedge_{j\neq k}(a_j\in u_j)\}$ is a subset of the set $\{x_{f(a_1,\ldots,a_n)}\midsp  \bigwedge_{j}(a_j\in u_j)\}$, hence taking joins it follows that $\widehat{f}(\vec{u}[w]_k)\subseteq\widehat{f}(\vec{u})$. By definition, if $xR\vec{u}$ holds, then we obtain $\widehat{f}(\vec{u}[w]_k)\subseteq\widehat{f}(\vec{u})\subseteq x$, hence $xR\vec{u}[w]_k$ holds as well. Similarly for the relation $S$.
\end{proof}
\end{lemma}

\begin{lemma}\rm\label{unified relational}
In the canonical frame, $xR\vec{u}$ holds iff $\forall\vec{a}\in L^n\;(\vec{a}\in\vec{u}\lra f(\vec{a})\in x)$. Similarly, $yS\vec{v}$ holds iff $\forall\vec{a}\in L^n\;(\vec{a}\in\vec{v}\lra h(\vec{a})\in y)$.
\end{lemma}
\begin{proof}
 By definition $xR\vec{u}$ holds iff $\widehat{f}(\vec{u})\subseteq x$, where $\widehat{f}(\vec{u})$, by its definition \eqref{canonical point operators defn} is the filter generated by the elements $f(\vec{a})$, for $\vec{a}\in\vec{u}$, hence clearly  $\vec{a}\in\vec{u}$ implies $f(\vec{a})\in x$. Similarly for the relation $S$.
\end{proof}

\begin{lemma}
\rm\label{same relation}
Where $R',S'$ are the Galois dual relations of the canonical relations $R,S$, $yR'\vec{u}$ holds iff $\widehat{f}(\vec{u})\upv y$ iff $\exists\vec{b}(\vec{b}\in\vec{u}\wedge f(\vec{b})\in y)$. Similarly, $xS'\vec{v}$ holds iff $x\upv \widehat{h}(\vec{v})$ iff $\exists\vec{e}(\vec{e}\in\vec{v}\wedge h(\vec{e})\in x)$.
\end{lemma}
\begin{proof}
  By definition of the Galois dual relation, $yR'\vec{u}$ holds iff for all $x\in X$, if $xR\vec{u}$ obtains, then $x\upv y$. By definition of the canonical relation $R$, for any $x\in X$, $xR\vec{u}$ holds iff $\widehat{f}(\vec{u})\subseteq x$ and thereby $\widehat{f}(\vec{u})R\vec{u}$ always obtains. Hence, $yR'\vec{u}$ is equivalent to $\forall x\in X\;(\widehat{f}(\vec{u})\subseteq x\lra x\cap y\neq\emptyset)$, from which it follows that $\widehat{f}(\vec{u})\upv y$ iff $yR'\vec{u}$ obtains.

To show that $yR'\vec{u}$ holds  iff $\exists\vec{a}(\vec{a}\in\vec{u}\wedge f(\vec{a})\in y)$, since the direction from right to left is trivially true, assume $yR'\vec{u}$, or, equivalently by the argument given above, assume that $\widehat{f}(\vec{u})\upv y$, i.e. $\widehat{f}(\vec{u})\cap y\neq \emptyset$ and let $e\in \widehat{f}(\vec{u})\cap y$. By $e\in\widehat{f}(\vec{u})$ and definition of $\widehat{f}(\vec{u})$ as the filter generated by the set $\{f(\vec{a})\midsp\vec{a}\in\vec{u}\}$, let $\vec{a}^1,\ldots,\vec{a}^s$, for some positive integer $s$, be $n$-tuples of lattice elements (where $\vec{a}^r=(a^r_1,\ldots,a^r_n)$, for $1\leq r\leq s$) such that $f(\vec{a}^1)\wedge\cdots\wedge f(\vec{a}^s)\leq e$ and $a^r_j\in u_j$ for each $1\leq r\leq s$ and $1\leq j\leq n$. Recall that the distribution type of $f$ is $\delta(f)=(i_1,\ldots,i_n;1)$, where for $j=1,\ldots,n$ we have $i_j\in\{1,\partial\}$ and define elements $b_1,\ldots,b_n$ as follows.
\[
b_j=\left\{
\begin{array}{cl}
a^1_j\wedge\cdots\wedge a^s_j &\;\mbox{ if }\; i_j=1\\[1mm]
a^1_j\vee\cdots\vee a^s_j &\;\mbox{ if }\; i_j=\partial
\end{array}
\right.
\]
When $i_j=1$, $f$ is monotone at the $j$-th argument place, $u_j$ is a filter and $b_j\leq a^r_j\in u_j$, for all $r=1,\ldots,s$, so that $b_j=a^1_j\wedge\cdots\wedge a^s_j\in u_j$. Similarly, when $i_{j'}=\partial$, $f$ is antitone at the $j'$-th argument place, while $u_{j'}$ is an ideal, so that $b_{j'}=a^1_{j'}\vee\cdots\vee a^s_{j'}\in u_{j'}$. This shows that $\vec{b}\in\vec{u}$ and it remains to show that $f(\vec{b})\in y$. We argue that $f(\vec{b})\leq f(\vec{a}^1)\wedge\cdots\wedge f(\vec{a}^s)\leq e$ and the desired conclusion follows by the fact that $e\in y$, an ideal.

For any $1\leq r\leq s$, let $\vec{a}^r[b_j]_j^{i_j=1}$ be the result of replacing $a^r_j$ by $b_j$ in the tuple $\vec{a}^r$ and in every position $j$ from 1 to $n$ such that $i_j=1$ in the distribution type of $f$. Since $b_j\leq a^r_j$ and $f$ is monotone at any such $j$-th argument place, it follows that $f(\vec{a}^r[b_j]_j^{i_j=1})\leq f(\vec{a}^r)$, for all $1\leq r\leq s$.

In addition, for any $1\leq r\leq s$, let $\vec{a}^r[b_j]_j^{i_j=1}[b_{j'}]_{j'}^{i_{j'}=\partial}$ be the result of replacing $a^r_{j'}$ by $b_{j'}$ in the tuple $\vec{a}^r[b_j]_j^{i_j=1}$ and in every position $j'$ from 1 to $n$ such that $i_{j'}=\partial$ in the distribution type of $f$. Since $b_{j'}\geq a^r_{j'}$ and $f$ is antitone at any such $j'$-th argument place, it follows that $f(\vec{a}^r[b_j]_j^{i_j=1}[b_{j'}]_{j'}^{i_{j'}=\partial})\leq f(\vec{a}^r[b_j]_j^{i_j=1})\leq f(\vec{a}^r)$, for all $1\leq r\leq s$. Since $\vec{a}^r[b_j]_j^{i_j=1}[b_{j'}]_{j'}^{i_{j'}=\partial}=\vec{b}$ we obtain  that
\[
f(\vec{b})=f(\vec{a}^r[b_j]_j^{i_j=1}[b_{j'}]_{j'}^{i_{j'}=\partial})\leq f(\vec{a}^r[b_j]_j^{i_j=1})\leq f(\vec{a}^1)\wedge\cdots\wedge f(\vec{a}^s)\leq e
\]
hence $f(\vec{b})\in y$ and  this completes the proof, as far as the relation $R$ is concerned.

The argument for the relation $S$ is similar and can be safely left to the reader.
\end{proof}

\begin{lemma}\rm
\label{canonical frame properties prop}
 In the canonical frame, all sections of the Galois dual relations $R', S'$ of the canonical relations $R,S$ are Galois sets.
\end{lemma}
\begin{proof}
There are two cases to handle, one for each  of the relations $R',S'$, with two subcases for each one, depending on whether $i_k$ is $1$, or $\partial$.\\

Case of the relation $R'$):  We have $R'\vec{u}=(R\vec{u})'$, by definition, so the section $R'\vec{u}$ is a Galois (co-stable) set. It remains to be shown that for any $y\in Y$ and $\vec{u}[\_]_k$, the $k$-th section $yR'\vec{u}[\_]_k$ is a Galois set, for any $1\leq k\leq n$. There are two subcases to consider, accordingly as $i_k=1$, or $i_k=\partial$ and recall that $\delta(f)=(i_1,\ldots,i_n;1)$. Hence if $i_k=1$, then $f$ is monotone and it distributes over finite joins at the $k$-th argument place and if $i_k=\partial$, then $f$ is antitone and it co-distributes at the $k$-th argument place over finite meets (turning them to joins).\\

Subcase $i_k=1$:\\ Then $yR'\vec{u}[\_]_k\subseteq X=\filt(\mathcal{L})$.

Let $v$ be the ideal generated by the set $V=\{b\in L\midsp \exists \vec{a}[\_]_k\;\; f(\vec{a}[b]_k)\in y\}$. For any $x\in X$ such that $yR'\vec{u}[x]_k$ holds, by Lemma \ref{same relation}  we have $\widehat{f}(\vec{u}[x]_k)\upv y$, equivalently, $\exists \vec{a}[\_]_k\exists b\;(\vec{a}[b]_k\in \vec{u}[x]_k\;\wedge\; f(\vec{a}[b]_k)\in y)$. Thus $b\in x\cap v$ and so $ yR'\vec{u}[\_]_k \upv v$.

We assume $z\in \left( yR'\vec{u}[\_]_k \right)^{\prime\prime}$ and we need to show that $yR'\vec{u}[z]_k$. The assumption implies that $z\upv v$, i.e. for some lattice element $e$ we have $e\in z\cap v$. By $e\in v$, let $b_1,\ldots,b_s\in V\subseteq v$, for some positive integer $s$, be elements such that $e\leq b_1\vee\cdots\vee b_s\in v$, since $v$ is an ideal.

Since $b_1,\ldots,b_s\in V$, there are tuples of lattice elements $\vec{c}^r[\_]_k$ such that $f(\vec{c}^r[b_r]_k)\in y$, for each  $1\leq r\leq s$. Considering the distribution type of $f$ and as in the proof of Lemma \ref{same relation}, define the tuple of elements $\vec{a}=(a_1,\ldots,a_n)$ by
\[
a_j=\left\{
\begin{array}{cl}
c^1_j\wedge\cdots\wedge c^s_j &\;\mbox{ if }\; i_j=1\\[1mm]
c^1_j\vee\cdots\vee c^s_j &\;\mbox{ if }\; i_j=\partial
\end{array}
\right.
\]
For each $1\leq r\leq s$ let $\vec{c}^r[b_r]_k[a_j]_j^{i_j=1}$ be the result of replacing in $\vec{c}^r[b_r]_k$ all $c^r_j$ by $a_j$ whenever $i_j=1$ in the distribution type of $f$. Then for each $r$ as above and by monotonicity of $f$ at the $j$-th argument place whenever $i_j=1$ we have $f(\vec{c}^r[b_r]_k[a_j]_j^{i_j=1})\leq
f(\vec{c}^r[b_r]_k)\in y$, an ideal, hence $f(\vec{c}^r[b_r]_k[a_j]_j^{i_j=1})\in y$. Let also $\vec{c}^r[b_r]_k[a_j]_j^{i_j=1}[a_{j'}]_{j'}^{i_{j'}=\partial}$ be the result of further replacing in  $\vec{c}^r[b_r]_k[a_j]_j^{i_j=1}$ all $c^r_{j'}$ by $a_{j'}$ whenever $i_{j'}=\partial$ in the distribution type of $f$. Since $f$ is antitone at the $j'$-th position when $i_{j'}=\partial$ (given that the output type of $f$ is assumed to be $i_{n+1}=1$), we obtain $f(\vec{c}^r[b_r]_k[a_j]_j^{i_j=1}[a_{j'}]_{j'}^{i_{j'}=\partial} )\leq f(\vec{c}^r[b_r]_k[a_j]_j^{i_j=1})\leq f(\vec{c}^r[b_r]_k)\in y$ and so $f(\vec{c}^r[b_r]_k[a_j]_j^{i_j=1}[a_{j'}]_{j'}^{i_{j'}=\partial} )\in y$. But, having performed substitutions in all places (except for the $k$-th) $\vec{c}^r[b_r]_k[a_j]_j^{i_j=1}[a_{j'}]_{j'}^{i_{j'}=\partial}=\vec{a}[b_r]_k$, for each $1\leq r\leq s$.

It follows from the above that, for each $1\leq r\leq s$ we have $f(\vec{a}[b_r]_k)\in y$, an ideal, hence also $\bigvee_{r=1}^{r=s}f(\vec{a}[b_r]_k)\in y$. By case assumption, $i_k=1$, hence $f$ distributes over finite joins in the $k$-th argument place and then we obtain that $f(\vec{a}[b_1\vee\cdots\vee b_s]_k)=\bigvee_{r=1}^{r=s}f(\vec{a}[b_r]_k)\in y$. By $e\leq b_1\vee\cdots\vee b_s$ and monotonicity of $f$ at the $k$-th argument place we obtain $f(\vec{a}[e]_k)\leq f(\vec{a}[b_1\vee\cdots\vee b_s]_k)\in y$, hence also $f(\vec{a}[e]_k)\in y$.

Therefore, there exists $\vec{a}[e]_k\in\vec{u}[z]_k$ such that $f(\vec{a}[e]_k)\in y$ which, by Lemma \ref{same relation}  is equivalent to $yR'\vec{u}[z]_k$. This shows that $(yR'\vec{u}[\_]_k)^{\prime\prime}\subseteq yR'\vec{u}[\_]_k$, so the section has been shown to be a Galois (stable) set.\\

The subcase $i_k=\partial$ and the case of the relation $S$ are treated similarly and we leave details to the interested reader.
\end{proof}

The canonical frame for a lattice expansion $\mathcal{L}=(L,\leq,\wedge,\vee,0,1,f,h)$, where $\delta(f)=(i_1,\ldots,i_n;1)$ and $\delta(h)=(t_1,\ldots,t_n;\partial)$ ($i_j,t_j\in\{1,\partial\}$) is  the structure $\mathcal{L}_+=\mathfrak{F}=(\filt(\mathcal{L}),\upv,\idl(\mathcal{L}),R,S)$. By Proposition \ref{canonical frame properties prop}, the canonical relations $R,S$ are compatible with the Galois connection generated by $\upv\;\subseteq X\times Y$, in the sense that all sections of their Galois dual relations are Galois sets.
Set operators $\alpha_R, \eta_S$ are defined as in Section \ref{image ops section} and we let $\overline{\alpha}_R,\overline{\eta}_S$ be the closures of their restrictions to Galois sets (according to their distribution types). Note that $\overline{\alpha}_R(\vec{F})\in\gpsi$, while $\overline{\eta}_S(\vec{G})\in\gphi$, given the output types of $f,h$ (alternatively, given the sort types of $R,S$).

It follows from Theorem \ref{distribution from stability thm} and Lemma \ref{canonical frame properties prop}, that the sorted operators $\overline{\alpha}_R,\overline{\eta}_S$ on Galois sets distribute over arbitrary joins of Galois sets (stable or co-stable, according to the sort types of $R,S$) in each argument place.

Note that $\overline{\alpha}_R,\overline{\eta}_S$ are sorted maps, taking their values in $\gpsi$ and $\gphi$, respectively. We define single-sorted maps on $\gpsi$ (analogously for $\gphi$) by composition with the Galois connection
\begin{eqnarray}
  \overline{\alpha}_f(A_1,\ldots,A_n) &=& \overline{\alpha}_R(\ldots,\underbrace{A_j}_{i_j=1},\ldots,\underbrace{A^\prime_r}_{i_r=\partial},\ldots) \hskip5mm(A_1,\ldots,A_n\in\gpsi)\label{single-sorted f} \\
  \overline{\eta}_h(B_1,\ldots,B_n) &=& \overline{\eta}_R(\ldots,\underbrace{B_r}_{i_r=\partial},\ldots,\underbrace{B^\prime_j}_{i_j=1},\ldots) \hskip5mm(B_1,\ldots,B_n\in\gphi)\label{single-sorted h}
\end{eqnarray}
Given that the Galois connection is a duality of Galois stable and Galois co-stable sets, it follows that the distribution type of $\overline{\alpha}_f$ is that of $f$ and that $\overline{\alpha}_f$ distributes, or co-distributes, over arbitrary joins and meets in each argument place, according to its distribution type, returning joins in $\gpsi$. Similarly, for $\overline{\eta}_h$.
Thus, the lattice representation maps $\zeta_1:(L,\leq,\wedge,\vee,0,1)\lra\gpsi$ and $\zeta_\partial:(L,\leq,\wedge,\vee,0,1)\lra\gphi$ are extended to maps $\zeta_1:\mathcal{L}\lra\gpsi$ and $\zeta_\partial:\mathcal{L}\lra \gphi$ by setting
\begin{eqnarray}
\zeta_1(f(a_1,\ldots,a_n))&=&\overline{\alpha}_f(\zeta_1(a_1),\ldots, \zeta_1(a_n))   
    \hskip2mm=\; \overline{\alpha}_R(\ldots,\underbrace{\zeta_1(a_j)}_{i_j=1},\ldots, \underbrace{\zeta_\partial(a_r)}_{i_r=\partial},\ldots)\nonumber\\
\zeta_\partial(f(a_1,\ldots,a_n)) &=& \left(\;\overline{\alpha}_f(\zeta_1(a_1),\ldots, \zeta_1(a_n))\;\right)^\prime
\label{canonical rep of lat ops f}\\[4mm]
\zeta_1(h(a_1,\ldots,a_n))&=& \left(\;\overline{\eta}_h(\zeta_\partial(a_1),\ldots,\zeta_\partial(a_n))\;\right)^\prime \nonumber\\
\zeta_\partial(h(a_1,\ldots,a_n))&=& \overline{\eta}_h(\zeta_\partial(a_1),\ldots,\zeta_\partial(a_n))\label{canonical rep of lat ops h}
\end{eqnarray}

It has been therefore established that there exists a map from normal lattice expansions to polarities with relations, as specified in the following concluding result.
\begin{coro}\rm
Given a normal lattice expansion $\mathcal{L}=(L,\leq,\wedge,\vee,0,1,\mathfrak{F}_1,\mathfrak{F}_\partial)$, where $\mathfrak{F}_1,\mathfrak{F}_\partial$ are families of normal lattice operators of output types 1 and $\partial$, respectively, the dual frame $\mathcal{L}_+$ of the lattice expansion $\mathcal{L}$ is a polarity with additional relations, where for a normal lattice operator $f$ of distribution type $(i_1,\ldots,i_n;i_{n+1})$ the corresponding frame relation $R_f$ is of sort type $(i_{n+1};i_1\cdots i_n)$ and where all sections of its Galois dual relation $R'_f$ are Galois sets.\telos
\end{coro}

\subsection{Representation, Canonical Extensions and RS-Frames}
\label{canext section}
A canonical lattice extension is defined in \cite{mai-harding} as a pair $(\alpha,C)$ where $C$ is a complete lattice and $\alpha$ is an embedding of a lattice $\mathcal{L}$ into $C$ such that the following density and compactness requirements are satisfied
\begin{quote}
$\bullet$ (density) $\alpha[{\mathcal L}]$ is {\em dense} in $C$, where the latter means that every element of $C$ can be expressed both as a meet of joins and as a join of meets of elements in $\alpha[{\mathcal L}]$\\[0.5mm]
$\bullet$ (compactness) for any set $A$ of closed elements and any set  $B$ of open elements of $C$, $\bigwedge A\leq\bigvee B$ iff there exist finite subcollections $A'\subseteq A, B'\subseteq B$ such that $\bigwedge A'\leq\bigvee B'$
\end{quote}
where the {\em closed elements} of $C$ are defined in \cite{mai-harding} as the elements in the meet-closure of the representation map $\alpha$ and the {\em open elements} of $C$ are defined dually as the join-closure of the image of $\alpha$.
\begin{prop}\rm
$\gpsi$ (the lattice of Galois stable subsets of the set of filters) is  a canonical extension of the (bounded) lattice $(L,\leq,\wedge,\vee,0,1)$.
\end{prop}
\begin{proof}
This was shown by Gehrke and Harding in \cite{mai-harding}. More precisely, existence of canonical extensions is proven in \cite{mai-harding} by demonstrating that the compactness and density requirements are satisfied in the representation due to Hartonas and Dunn \cite{sdl}, which is precisely the representation presented in Section \ref{canonical frame section}.
\end{proof}

\begin{prop}\rm
The canonical representations of the normal lattice operators $f, h$, of respective output types $1,\partial$, as defined by the equations \eqref{canonical rep of lat ops f} and \eqref{canonical rep of lat ops h}, are the $\sigma$ and $\pi$-extension, respectively, of $f,h$.
\end{prop}
\begin{proof}
In the representation of Section \ref{canonical frame section}, the closed and open elements of $\gpsi$ are  the sets of the form $\Gamma x\; (x\in X)$ and ${}\rperp\{y\}\;(y\in Y)$, respectively.  For a unary lattice operator $f:\mathcal{L}\lra\mathcal{L}$,  its $\sigma$-extension in a canonical extension $C$ of the lattice $\mathcal{L}$ is defined in \cite{mai-harding} by equation \eqref{sigma-extn}, where $\type{K}$ is the set of closed elements of $C$ and $\type{O}$ is its set of open elements
\begin{eqnarray}
f_\sigma(k)=\bigwedge\{f(a)\midsp k\leq a\in L\} &\hskip1cm& f_\sigma(u)=\bigvee\{f_\sigma(k)\midsp\type{K}\ni k\leq u\}\label{sigma-extn}\\
f_\pi(o)=\bigvee \{f(a)\midsp L\ni a\leq o\} && f_\pi(u)=\bigwedge\{f_\pi(o)\midsp u\leq o\in\type{O}\}\label{pi-extn}
\end{eqnarray}
where in these definitions $\cal L$ is identified with its isomorphic image in $C$ and $a\in {\cal L}$ is then identified with its representation image.

Working concretely with the canonical extension of \cite{sdl}, the $\sigma$ extension $f_\sigma:{\mathcal L}_\sigma\lra{\mathcal L}_\sigma$ of a monotone map $f$ as in equation (\ref{sigma-extn}) and the dual $\sigma$-extension $f_\sigma^\partial:{\mathcal L}^\partial_\sigma\lra{\mathcal L}^\partial_\sigma$ (not used in \cite{mai-harding}) are defined by instantiating equation (\ref{sigma-extn}) in the concrete canonical extension of \cite{sdl} considered here by setting, for $x\in X$ and $y\in Y$ and where $x_e$ is a principal filter, $y_e$ a principal ideal and the closed elements are precisely the principal upper sets $\Gamma u$ (for each of $X$, $Y$).
\begin{tabbing}
$f_\sigma(\Gamma x)$\hskip2mm\==\hskip2mm\=$\bigwedge\{\alpha_X(fa)\midsp a\in{\mathcal L}, \Gamma x\leq\alpha_X(a)\}$\hskip2mm\==\hskip2mm\=$\bigwedge\{\Gamma x_{fa}\midsp\Gamma x\subseteq\Gamma x_a\}$\\
\>=\>$\bigwedge\{\Gamma x_{fa}\midsp a\in x\}$\>=\>$\Gamma(\bigvee\{x_{fa}\midsp a\in x\})$\\
\>=\> $\Gamma(\widehat{f}(x))$\\[1mm]
$f_\sigma^\partial(\Gamma y)$\>=\>$\bigwedge\{\alpha_Y(fa)\midsp a\in{\mathcal L}, \Gamma y\leq\alpha_Y(a)\}$\>=\>$\bigwedge\{\Gamma y_{fa}\midsp\Gamma y\subseteq\Gamma y_a\}$\\
\>=\>$\bigwedge\{\Gamma y_{fa}\midsp a\in y\}$\>=\>$\Gamma(\bigvee\{y_{fa}\midsp a\in y\})$
\end{tabbing}
Hence $f_\sigma(\Gamma x)=\Gamma(\widehat{f}(x))$, where $\widehat{f}$ is the point operator we defined, after \cite{dloa}, by equation \eqref{canonical point operators defn}. For an $n$-ary operator $f$, first observe that in a product $\prod_{j=1}^{j=n}\mathcal{G}(Z_{i_j})$, where when $i_j=\partial$ then $Z_{i_j}=Y$ and $\gphi$, closed elements are $n$-tuples of closed elements of the factors $(\Gamma u_1,\ldots,\Gamma u_n)$. Then $f_\sigma(\Gamma u_1,\ldots,\Gamma u_n)=\Gamma(\widehat{f}(\vec{u}))$, by the same analysis. This is a sorted operator and by composition with the Galois connection we obtain (in the case examined the output type is 1) the single-sorted $\sigma$-extension $f^\sigma(\ldots,\underbrace{\Gamma u_j}_{i_j=1},\ldots,\underbrace{{}\rperp\{u_r\}}_{i_r=\partial},\ldots)=\Gamma(\widehat{f}(\vec{u}))$.

For an arbitrary Galois stable set $A$ and unary monotone $f$, $f^\sigma(A)$ is defined in \eqref{sigma-extn} using join-density of closed elements. Hence we obtain $f^\sigma(A)=\bigvee_{x\in A}f_\sigma(\Gamma x)=\bigvee_{x\in A}\Gamma(\widehat{f}(x))$.

For an $n$-ary monotone map we similarly obtain that $f^\sigma(\vec{F})=\bigvee_{\vec{u}\in\vec{F}}\Gamma(\widehat{f}(\vec{u}))$. Since $w\in \Gamma(\widehat{f}(\vec{u}))$ iff $\widehat{f}(\vec{u})\subseteq w$ iff $wR\vec{u}$, by the way the canonical relation $R$ was defined in equation \eqref{canonical relations defn}, so that $\Gamma(\widehat{f}(\vec{u}))=R\vec{u}$, we obtain that (observing that $R\vec{u}$ is a Galois set, indeed a closed element of $\gpsi$)
\[
f^\sigma(\vec{F})=\bigvee_{\vec{u}\in\vec{F}}R\vec{u}=\left(\bigcup_{\vec{u}\in\vec{F}}R\vec{u}\right)^{\prime\prime}= \overline{\alpha}_R(\vec{F})
\]
Note that $\sigma$-extensions, as defined in \cite{mai-harding}, are sorted maps and then a single-sorted map is obtained by composing with the Galois connection, as shown in equations \eqref{canonical rep of lat ops f} and \eqref{canonical rep of lat ops h}.

The $\pi$-extension is simply the Galois image of the dual $\sigma$-extension, so there is nothing new to discuss and the proof is complete. We only note further that the way we have canonically proceeded is to represent a lattice operator with output type 1 by its $\sigma$-extension and one of output type $\partial$ by its $\pi$-extension.
\end{proof}

\begin{rem}[Canonical Relations in the RS-Frames Approach]\rm
In modeling the Lambek calculus product operator $\circ$, of distribution type $(1,1;1)$ (see Example \ref{example and strategy ex 3}), the canonical relation $R^{111}$ was defined by (using the point operators) $uRxz$ iff $x\overt z\subseteq u$, where $x\overt z$ is the filter generated by the elements $a\circ c$, with $a\in x$ and $c\in z$. By Lemma \ref{unified relational}, specialized to this case, this amounts to the classical definition of a canonical relation, familiar from the Boolean and distributive case, by the clause
\[
uRxz\;\mbox{ iff }\; \forall a,c\in\mathcal{L}\;(a\in x\;\wedge\;c\in z\;\lra\; a\circ c\in u)
\]
By Lemma \ref{same relation}, specialized to the particular case, $yR'xz$ holds iff $x\overt z\upv y$ iff $\exists a,c\in\mathcal{L}(a\in x\;\wedge\;c\in z\;\wedge a\circ c\in y)$, which is precisely
Gehrke's \cite{mai-gen} canonical relation definition for the Lambek product operator.

This is no isolated matter, as in the RS-frames approach the choice is made to work directly with a relation that can be nevertheless defined as the Galois dual of a classically defined accessibility relation. This is also witnessed by the way Goldblatt \cite{goldblatt-morphisms2019} proceeds, generalizing on Gehrke's \cite{mai-gen}, to define relations and set-operators in a frame. Indeed, examining  (for his case of interest) additive operators $F$ on stable sets he defines a relation $S_F$ by setting $yS_F\vec{z}$ iff $F(\Gamma z_1,\ldots, \Gamma z_n)\upv y$ (which in the case of a binary, completely additive operator $F$, is equivalent in the canonical frame to $z_1\overt z_2\upv y$). It is merely a matter of choice and convenience, given the purpose at hand, which relation to decide to work with. Gehrke \cite{mai-gen} does indeed point out that instead of using the relation $S\subseteq Y\times(X\times X)$ defined as above (for the Lambek product operator), one could use a relation $R\subseteq X\times(X\times X)$, which is actually the Galois dual of $S$, but she does not dwell much on the matter.

Though, to the best of this author's knowledge, it has not been made explicit in the RS-frames approach how relations are to be defined corresponding to arbitrary normal lattice operators in general  (but only in cases of particular examples), the relations on an RS-frame corresponding to normal lattice operators of some distribution type $(i_1,\ldots,i_n;i_{n+1})$ are the Galois duals of our canonical accessibility relations, hence they are systematically of sort type $(\overline{i_{n+1}};i_1\cdots i_n)$, and operators are defined from them. For example, Goldblatt \cite{goldblatt-morphisms2019} (generalizing Gehrke's \cite{mai-gen} definition for the Lambek product operator) defines from a relation $S\subseteq Y\times X^n$ an operator $\type{F}_S$ by setting
\begin{eqnarray}
    \type{F}_S^\bullet(\vec{F}) &=& \bigcap\{S\vec{z}\midsp \vec{z}\in\vec{F}\} \\
    \type{F}_S(\vec{F}) &=& {}\rperp \type{F}_S^\bullet(\vec{F})\;=\; \bigvee\{{}\rperp(S\vec{z})\midsp \vec{z}\in\vec{F})\}=\bigvee_{\vec{z}\in\vec{F}}{}\rperp(S\vec{z})\label{Goldblatt op}
  \end{eqnarray}
The relation $R$ defined as the Galois dual of $S$, i.e. by $R\vec{z}={}\rperp(S\vec{z})$ is precisely a relation of sort type $(i_{n+1};i_1\cdots i_n)$ and, assuming section stability, $R$ and $S$ are each other's Galois dual. Therefore we obtain
\begin{eqnarray}
\type{F}_S(\vec{F}) &=&\bigvee_{\vec{z}\in\vec{F}}{}\rperp(S\vec{z})\;=\;\bigvee_{\vec{z}\in\vec{F}}R\vec{z}\;=\; \left(\bigcup_{\vec{z}\in\vec{F}}R\vec{z}\right)^{\prime\prime}\label{Goldblatt op 2}
\end{eqnarray}

A comparison of equations \eqref{closure of restriction} and \eqref{Goldblatt op 2} reveals then that the two definitions are variants of each other.
\end{rem}

\section{Conclusions}
We have argued that the two approaches, the one developed in this article (concluding and completing our previous recent work on the subject) and the RS-frames approach really only differ in whether the polarity is assumed to be separated and reduced or not. The results of this article have shown that nothing is lost by dropping these additional assumptions, as far as the semantics of logics without distribution is concerned.

There are three points of interest, however, that are worth making.

First, a Stone type duality for RS-frames (essentially for Hartung's lattice representation) has encountered difficulties, similar to these encountered in extending Urquhart's representation to a full Stone duality. In \cite{sdl-exp}, we have developed a duality result for normal lattice expansions, extending the representation of \cite{sdl}. In view of Goldblatt's recent proposal \cite{goldblatt-morphisms2019} of a notion of bounded morphisms for polarities, this result, combined with the results of this article, can be improved to a Stone duality for normal lattice expansions with bounded morphisms as the morphisms in the dual category of polarities. We shall not dwell on this more here, for lack of space.

Second, the approach we have presented in this article allows for relating the logic of non-distributive lattices to the sorted, residuated (poly)modal logic of polarities with relations, where the residuated pair of modal operators is generated by the complement of the Galois relation of the frame. Preliminary results in this direction have been reported in \cite{pll7,redm} by this author, but the area is far from fully explored. Regarding non-distributive logics as fragments of sorted, residuated (poly)modal logics allows for importing techniques and results from modal logic in the field of logics lacking distribution.

Finally, we believe that the semantic framework presented in this article fully complements Dunn's gaggle theory project and in an important sense it completes the project for the case of non-distributive logical calculi.

\bibliographystyle{plain}

\begin{thebibliography}{10}

\bibitem{dunnLL}
Gerard Allwein and J.~Michael Dunn.
\newblock Kripke models for linear logic.
\newblock {\em The Journal of Symbolic Logic}, 58(2):514--545, 1993.

\bibitem{gglbook}
Katalin Bimb\'{o} and J.~Michael Dunn.
\newblock {\em Generalized Galois Logics. Relational Semantics of Nonclassical
  Logical Calculi}, volume 188.
\newblock CSLI Lecture Notes, CSLI, Stanford, CA, 2008.

\bibitem{birkhoff}
Garrett Birkhoff.
\newblock {\em Lattice theory}.
\newblock American Mathematical Society Colloquium Publications 25, American
  Mathematical Society, Providence, Rhode Island, third edition, 1979.
\newblock (corrected reprint of the 1967 third edition).

\bibitem{mai-grishin}
Anna Chernilovskaya, Mai Gehrke, and Lorijn van Rooijen.
\newblock Generalised {K}ripke semantics for the {L}ambek-{G}rishin calculus.
\newblock {\em Logic Journal of the IGPL}, 20(6):1110--1132, 2012.

\bibitem{palmigiano-categories}
Willem Conradie, Sabine Frittella, Alessandra Palmigiano, Michele Piazzai,
  Apostolos Tzimoulis, and Nachoem Wijnberg.
\newblock Categories: How {I} learned to stop worrying and love two sorts.
\newblock In {\em Logic, Language, Information, and Computation - 23rd
  International Workshop, WoLLIC 2016, Puebla, Mexico, August 16-19th, 2016.
  Proceedings}, pages 145--164, 2016.

\bibitem{conradie-palmigiano}
Willem Conradie and Alessandra Palmigiano.
\newblock Algorithmic correspondence and canonicity for non-distributive
  logics.
\newblock {\em Ann. Pure Appl. Logic}, 170(9):923--974, 2019.

\bibitem{conradie2021nondistributive}
Willem Conradie, Alessandra Palmigiano, Claudette Robinson, and Nachoem
  Wijnberg.
\newblock Non-distributive logics: from semantics to meaning, 2021.

\bibitem{conradie2018goldblattthomason}
Willem Conradie, Alessandra Palmigiano, and Apostolos Tzimoulis.
\newblock Goldblatt-{T}homason for {LE}-logics, 2018.

\bibitem{COUMANS201450}
Dion Coumans, Mai Gehrke, and Lorijn van Rooijen.
\newblock Relational semantics for full linear logic.
\newblock {\em Journal of Applied Logic}, 12(1):50 -- 66, 2014.
\newblock Logic Categories Semantics.

\bibitem{dunn-ggl}
J.~Michael Dunn.
\newblock Gaggle theory: An abstraction of {G}alois coonections and resuduation
  with applications to negations and various logical operations.
\newblock In {\em Logics in AI, Proceedings of European Workshop JELIA 1990,
  LNCS 478}, pages 31--51, 1990.

\bibitem{dunn-gehrke}
Jon~Michael Dunn, Mai Gehrke, and Alessandra Palmigiano.
\newblock Canonical extensions and relational completeness of some
  substructural logics.
\newblock {\em Journal of Symbolic Logic}, 70:713--740, 2005.

\bibitem{vac-et-al}
Ivo D\"{u}ntsch, Ewa Orlowska, Anna~Maria Radzikowska, and Dimiter Vakarelov.
\newblock Relational representation theorems for some lattice-based structures.
\newblock {\em Journal of Relational Methods in Computer Science (JORMICS)},
  1:132--160, 2004.

\bibitem{wille2}
Bernhard Ganter and Rudolph Wille.
\newblock {\em Formal Concept Analysis: Mathematical Foundations}.
\newblock Springer, 1999.

\bibitem{mai-gen}
Mai Gehrke.
\newblock Generalized {K}ripke frames.
\newblock {\em Studia Logica}, 84(2):241--275, 2006.

\bibitem{mai-harding}
Mai Gehrke and John Harding.
\newblock Bounded lattice expansions.
\newblock {\em Journal of Algebra}, 238:345--371, 2001.

\bibitem{mai-jons}
Mai Gehrke and Bjarni J\'{o}nsson.
\newblock Bounded distributive lattice expansions.
\newblock {\em Math. Scand.}, 94(2):13--45, 2004.

\bibitem{goldb}
Robert Goldblatt.
\newblock Semantic analysis of orthologic.
\newblock {\em Journal of Philosophical Logic}, 3:19--35, 1974.

\bibitem{goldblatt-morphisms2019}
Robert Goldblatt.
\newblock Morphisims and duality for polarities and lattices with operators.
\newblock {\em {FLAP}}, 7(6):1017--1070, 2020.

\bibitem{grishin}
V.N. Grishin.
\newblock On a generalization of the {A}jdukiewicz-{L}ambek system.
\newblock In A.I. Mikhailov, editor, {\em Studies in Nonclassical Logics and
  Formal Systems}, pages 315--334. Nauka, 1983.

\bibitem{dloa}
Chrysafis Hartonas.
\newblock Duality for lattice-ordered algebras and for normal algebraizable
  logics.
\newblock {\em Studia Logica}, 58:403--450, 1997.

\bibitem{sdl-exp}
Chrysafis Hartonas.
\newblock Stone {D}uality for {L}attice {E}xpansions.
\newblock {\em Oxford Logic Journal of the {IGPL}}, 26(5):475--504, 2018.

\bibitem{discr}
Chrysafis Hartonas.
\newblock Discrete duality for lattices with modal operators.
\newblock {\em J. Log. Comput.}, 29(1):71--89, 2019.

\bibitem{discres}
Chrysafis Hartonas.
\newblock Duality results for (co)residuated lattices.
\newblock {\em Logica Universalis}, 13(1):77--99, 2019.

\bibitem{gts}
Chrysafis Hartonas.
\newblock {Game-theoretic semantics for non-distributive logics}.
\newblock {\em Logic Journal of the IGPL}, 27(5):718--742, 01 2019.

\bibitem{pll7}
Chrysafis Hartonas.
\newblock Lattice logic as a fragment of (2-sorted) residuated modal logic.
\newblock {\em Journal of Applied Non-Classical Logics}, 29(2):152--170, 2019.

\bibitem{redm}
Chrysafis Hartonas.
\newblock Modal translation of substructural logics.
\newblock {\em Journal of Applied Non-Classical Logics}, 30(1):16--49, 2020.

\bibitem{iulg}
Chrysafis Hartonas and J.~Michael Dunn.
\newblock Duality theorems for partial orders, semilattices, galois connections
  and lattices.
\newblock Technical Report IULG-93-26, Indiana University Logic Group, 1993.

\bibitem{sdl}
Chrysafis Hartonas and J.~Michael Dunn.
\newblock Stone duality for lattices.
\newblock {\em Algebra Universalis}, 37:391--401, 1997.

\bibitem{pnsds}
Chrysafis Hartonas and Ewa Or{\l}owska.
\newblock Representation of lattices with modal operators in two-sorted frames.
\newblock {\em Fundam. Inform.}, 166(1):29--56, 2019.

\bibitem{hartung}
Gerd Hartung.
\newblock A topological representation for lattices.
\newblock {\em Algebra Universalis}, 29:273--299, 1992.

\bibitem{jt1}
Bjarni J\'{o}nsson and Alfred Tarski.
\newblock Boolean algebras with operators {I}.
\newblock {\em American Journal of Mathematics}, 73:891--939, 1951.

\bibitem{jt2}
Bjarni J\'{o}nsson and Alfred Tarski.
\newblock Boolean algebras with operators {II}.
\newblock {\em American Journal of Mathematics}, 74:8127--162, 1952.

\bibitem{residBA}
Bjarni J\'{o}nsson and Constantine Tsinakis.
\newblock Relation algebras as residuated boolean algebras.
\newblock {\em Algebra Universalis}, 30:469--478, 1993.

\bibitem{ewa-rfca}
Ewa Or{\l}owska.
\newblock Modal logic in the theory of information systems.
\newblock {\em Mathematical Logic Quarterly}, 30(13‐16):213--222.

\bibitem{hilary}
Hilary Priestley.
\newblock Representation of distributive lattices by means of ordered {S}tone
  spaces.
\newblock {\em Bull. Lond. Math. Soc.}, 2:186--190, 1970.

\bibitem{stone2}
Marshall~Harvey Stone.
\newblock Topological representation of distributive lattices and brouwerian
  logics.
\newblock {\em Casopsis pro Pestovani Matematiky a Fysiky}, 67:1--25, 1937.

\bibitem{stone1}
Marshall~Harvey Stone.
\newblock The representation of boolean algebras.
\newblock {\em Bull. Amer. Math. Soc.}, 44:807--816, 1938.

\bibitem{Suzuki-polarity-frames}
Tomoyuki Suzuki.
\newblock On polarity frames: Applications to substructural and lattice-based
  logics.
\newblock In Rajeev Gor\'e, Barteld Kooi, and Agi Kurucz, editors, {\em
  Advances in Modal Logic, Volume 10}, pages 533--552. CSLI Publications, 2014.

\bibitem{urq}
Alasdair Urquhart.
\newblock A topological representation of lattices.
\newblock {\em Algebra Universalis}, 8:45--58, 1978.

\bibitem{rough-vakarelov}
Dimiter Vakarelov.
\newblock {\em Information Systems, Similarity Relations and Modal Logics},
  chapter~16, pages 492--550.
\newblock Studies in Fuzziness and Soft Computing. Heidelberg: Physica-Verlag,
  1998.

\end{thebibliography}

\end{document}